\documentclass[12pt]{amsart}
\usepackage{amsfonts,amsmath,amssymb,amsthm,latexsym,verbatim}
\input xypic
\newtheorem{theorem}{Theorem}[section]

\theoremstyle{remark}

\theoremstyle{definition}

\DeclareMathOperator{\rad}{{\mathrm{rad}}}

\newcommand{\dprod}[2]{\langle {#1} , {#2}\rangle}	

\DeclareMathOperator{\Ch}{{\mathrm {Ch}}} 	
 	
\begin{document}
\title[Weyl modules for $E_6$]{Some Weyl modules of the algebraic groups of type $E_6$}
\author{Peter Sin}
\address{Department of Mathematics\\University of Florida\\ P. O. Box 118105\\ Gainesville
FL 32611\\ USA}
\date{}
\thanks{This work was partially supported by a grant from the Simons Foundation (\#204181 to Peter Sin)}


\begin{abstract}
Let $G$ be a simple algebraic group of type $E_6$ over an algebraically closed field
of characteristic $p>0$. We determine the submodule structure of the Weyl modules
with highest weight $r\omega_1$ for $0\leq r\leq p-1$, where 
$\omega_1$ is the fundamental weight of the standard $27$-dimensional
module. In the process, the structures of other Weyl modules with
highest weights linked to $r\omega_1$ are also found. 
\end{abstract}

\maketitle
\section{Introduction}
In this note we study certain Weyl modules for a simple, simply connected 
algebraic group $G$ of type $E_6$ over an algebraically closed field of characteristic $p>0$.
The modules we consider are for highest weights which are of the form $r\omega_1$, $0\leq r\leq p-1$, where $\omega_1$ is the highest weight of the ``standard'' $27$-dimensional
module, and we will give a full description of their $G$-submodules. If $P$ is the maximal parabolic subgroup stabilizing the highest weight vector in the 27-dimensional module $H^0(\omega_1)^*$, then the embedding of the projective variety $G/P$ for the associated line bundle is projectively normal \cite{RR}, so the homogeneous coordinate ring is
$
\bigoplus_{r\ge 0}H^0(r\omega_1).
$

As a consequence of Steinberg's Tensor Product Theorem \cite{St},  our results
also describe the simple $G$-socles of the modules $H^0(r\omega_1)$ for all $r\geq0$.

Our labelling of the fundamental roots and weights is according to Figure~\ref{dynk}.
\begin{figure}[ht]\label{dynk}
\centering
\begin{picture}(100,100)(50,-50)
\thicklines
   \put(5,30){$\alpha_1$}
   \put(55,30){$\alpha_2$}
   \put(105,30){$\alpha_3$}
   \put(155,30){$\alpha_5$}
   \put(205,30){$\alpha_6$}
   \put(105,-45){$\alpha_4$}
   \put(10,20){\line(1,0){200}}
  \put(110,20){\line(0,-1){50}}
  \multiput(10,20)(50,0){5}{\circle*{5}}
  \put(110,-30){\circle*{5}}
\end{picture}
\caption{}
\end{figure}
We describe the $E_6$ root system as follows. Let $e_i$, $i=1,\ldots,8$
be an orthonomal basis of an $8$-dimensional Euclidean space.
Then, in coordinates,  our root system $R$ is the union of the set
$$
\{ \pm e_i\pm e_j\mid 4\leq i<j\leq 8\}
$$
with the set
$$
\{\pm\frac12[(e_1-e_2-e_3)+\sum_{i=5}^8\pm e_i]\mid \text{number of minus signs is even}\}.
$$
A set of fundamental roots is
\begin{multline*}
S=\{\alpha_1=e_4-e_5, \alpha_2=e_5-e_6, \alpha_3=e_6-e_7, \alpha_4=e_7+e_8,\\
 \alpha_5=e_7-e_8,\alpha_6=\frac12(e_1-e_2-e_3-e_4-e_5-e_6-e_7+e_8)\}.
\end{multline*}

The fundamental dominant weights have coordinates
$$
\begin{aligned}
\omega_1&=\frac13(1,-1,-1,3,0,0,0,0), \quad
\omega_2=\frac13(2,-2,-2,3,3,0,0,0),\quad
\omega_3=(1,-1,-1,1,1,1,0,0),\\
\omega_4&=\frac12(1,-1,-1,1,1,1,1,1),\quad
\omega_5=\frac16(5,-5,-5,3,3,3,3,-3),\quad
\omega_6=\frac13(2,-2,-2,0,0,0,0,0).
\end{aligned}
$$ 

Our notation will be standard, following
\cite{Jantzen}. In particular we denote the Weyl module with highest weight
$\lambda$ by $V(\lambda)$ and its simple quotient by $L(\lambda)$. 
By definition, $V(\lambda)=H^0(-w_0\lambda)^*$, where $w_0$ is the longest element of the Weyl group \cite[II.2.13]{Jantzen}. Also,  $V(\lambda)\cong\,   ^\tau{H^0(\lambda)}$,
for a certain  anti-automorphism $\tau$ of $G$ that induces
the identity map on characters \cite[II. 2.12]{Jantzen}.
As $-w_0\omega_1=\omega_6$, the submodule structure of $V(r\omega_1)$ will yield the submodule structures of $V(r\omega_6)$,  $H^0(r\omega_1)$ and $H^0(r\omega_6)$ by applying $\tau$ and duality.

\begin{theorem}\label{main} Let $G$ be a simply connected, semisimple algebraic
group of type $E_6$ over an algebraically closed field of characteristic $p$.
The following statements give a complete description
of the submodule structure of the module $V(r\omega_1)$, for $0\leq r\leq p-1$. 
\begin{enumerate} 
\item[(a)] For $0\leq r\leq p-4$ the Weyl  module $V(r\omega_1)$ is simple.
\item[(b)]($r=p-3$)
\begin{enumerate}
\item[(i)]  If $p=3$, then $V((p-3)\omega_1)=V(0)$ is simple.
\item[(ii)] If $p=5$, there is an exact sequence
$$
0\to V(\omega_6)\to V(2\omega_1)\to L(2\omega_1)\to 0.
$$
\item[(iii)]If $p=7$, there is  an exact sequence
$$
0\to V(\omega_1+\omega_4)\to V(2\omega_1+\omega_6)\to V(4\omega_1)\to L(4\omega_1)\to 0.
$$
\item[(iv)]For $p\geq 11$, there is an exact sequence
\begin{multline*} 
0\to V((p-9)\omega_1)\to V((p-8)\omega_1+\omega_6)\\
\to V((p-8)\omega_1+\omega_2)\to V((p-6)\omega_1+\omega_4)\\
\to V((p-5)\omega_1+\omega_6)\to V((p-3)\omega_1)\to L((p-3)\omega_1)\to 0
\end{multline*}
\end{enumerate}
\item[(c)] ($r=p-2$) 
\begin{enumerate}
\item[(i)] If $p=2$ or $p=3$ the Weyl  module $V((p-2)\omega_1)$ is simple.
\item[(ii)] If $p=5$ there is an exact sequence
$$
0\to V(0)\to V(3\omega_1)\to L(3\omega_1)\to 0.
$$
\item[(iii)] If $p=7$, there is an exact sequence
$$
0\to V(\omega_4+\omega_6)\to V(\omega_1+2\omega_6)\to V(5\omega_1)\to L(5\omega_1)\to 0.
$$
\item[(iv)] For $p\geq 11$ there is an exact sequence
\begin{multline*} 
0\to V((p-10)\omega_1+\omega_2)\to V((p-9)\omega_1+\omega_5)\\
\to V((p-8)\omega_1+\omega_3)\to V((p-7)\omega_1+\omega_4+\omega_6)\\
\to V((p-6)\omega_1+2\omega_6)\to V((p-2)\omega_1)\to L((p-2)\omega_1)\to 0
\end{multline*}
\end{enumerate}
\item[(d)]($r=p-1$)
\begin{enumerate}
\item[(i)] If $p\leq 5$ the  the Weyl  module $V((p-1)\omega_1)$ is simple.
\item[(ii)] If $p=7$, there is an exact sequence
\begin{equation*}
0\to V(3\omega_6)\to V(6\omega_1)\to L(6\omega_1)\to 0.
\end{equation*}
\item[(iii)] For $p\geq 11$ there is an exact sequence
\begin{multline*} 
0\to V((p-11)\omega_1+2\omega_2)\to V((p-10)\omega_1+\omega_2+\omega_5)\\
\to V((p-9)\omega_1+\omega_3+\omega_6)\to V((p-8)\omega_1+\omega_4+2\omega_6)\\
\to V((p-7)\omega_1+3\omega_6)\to V((p-1)\omega_1)\to L((p-1)\omega_1)\to 0
\end{multline*}
\end{enumerate}
\item[(e)] In each of the above sequences the first and last nonzero terms
are simple modules and the other terms have two composition factors.
\end{enumerate}
\end{theorem}

We shall apply the {\it Jantzen Sum Formula} \cite[II.8.19]{Jantzen}
\begin{footnote}
{The validity of the sum formula for all $p$ was proved by Andersen.}
\end{footnote}:
The Weyl module $V(\lambda)$ has a descending filtration,
of submodules $V(\lambda)^i$, $i>0$, such that
\begin{equation*}
V(\lambda)^1=\rad(V(\lambda)),\quad\text{(so that $V(\lambda)/V(\lambda)^1\cong L(\lambda)$)}
\end{equation*}
and 
\begin{equation*}\label{jsum}
J(\lambda):=\sum_{i>0}\Ch(V(\lambda)^i)= -\sum_{\alpha>0}
\sum_{\{m: 0<mp<\dprod{\lambda+\rho}{\alpha^\vee}\}}v_p(mp)\chi(\lambda-mp\alpha)
\end{equation*}
We shall refer to the quantity $J(\lambda)$ as the Jantzen sum for $\lambda$ (or
for $V(\lambda)$). 
We recall that the weight $\rho$ is the half-sum of the positive roots
and $v_p(m)$ denotes the exponent of $p$ in the prime factorization of $m$.
Finally, the formal character $\chi(\mu)$ is the so-called Weyl character, defined 
in \cite[II.5.7]{Jantzen}, which has the following concrete description. 
There is a unique dominant weight of the form $w(\mu+\rho)$ , where $w\in W$. Let
$\mu'=w(\mu+\rho)-\rho$. Then $\chi(\mu)$ is equal to $\text{sign}(w)\Ch(V(\mu'))$ if $\mu'$ is dominant, and zero otherwise.  In particular $\chi(\mu)=\Ch(V(\mu))$ if $\mu$ is dominant and $\chi(\mu)=0$ if and only if $\lambda+\rho-mp\alpha$ is orthogonal to some root. 

To aid our computation, when $p$ and $\lambda\in X_+$ have been fixed,
we shall say that a multiple $m\alpha$ of a positive root is {\it relevant} if 
$0<mp<\dprod{\lambda+\rho}{\alpha^\vee}\}$ and that $m\alpha$ is a {\it contributor} 
if $m\alpha$ is relevant and $\chi(\lambda-mp\alpha)\neq 0$. 
We will call the quantity $v_p(mp)\chi(\lambda-mp\alpha)=\text{sign}(w)\chi(\mu')$ a {\it contribution} and the dominant weight $w(\lambda+\rho-mp\alpha)-\rho$ a {\it contributing weight}.
  
Thus, in computing the Jantzen sums
we can begin by determining the relevant root multiples, then determine which among them
is a contributor and finally add up the contributions.

We note that if $\alpha_0$ is the highest root, then $\langle r\omega_1+\rho, \alpha_0^\vee\rangle=r+11<2p$, when $p\geq 11$, so the only relevant root multiples are
actually roots. For the primes $p=2$, $3$, $5$ and $7$, we have to
take into account higher multiples.  

\subsection{Discussion of the proof of Theorem~\ref{main}}
The proof is by means of computations, whose results are compiled
in the tables below. The tables all have the same form. In the first column
are dominant weights $\lambda$. In the second column, for each $\lambda$ 
we list all the relevant root multiples. A root multiple is given by the tuple of coefficients in its expression as a sum of fundamental roots.
Since we are dealing throughout with a single root system $E_6$, the
relevant root multiples for any given weight are easily computed. 
The relevant root multiples for each weight $\lambda$ are divided
into non-contributors and contributors. 
For those root multiples $m\alpha$ which we claim to be noncontributors, 
we must exhibit a (co)root $\beta$ that is orthogonal to $\lambda+\rho-pm\alpha$. Note that $\beta$ is not necessarily unique up to a sign. The third column gives the coordinate tuple of the weight $\lambda+\rho-pm\alpha$ with respect to the fundamental weights and the fourth columns gives the coordinate tuple of a suitable $\beta$  with respect to the fundamental (co)roots. The reader can immediately verify that $\lambda+\rho-pm\alpha$ and $\beta$ are orthogonal, hence that $m\alpha$ is indeed a non-contributor.

For a contributing root multiple $m\alpha$, the third
column has the Weyl group element $w$ such that $w(\lambda+\rho-pm\alpha)$ 
is dominant and fourth column has the contributing weight $w(\lambda+\rho-pm\alpha)-\rho$.
The element $w$ is given as a tuple of indices $[i_1,i_2,\ldots, i_r]$, where
$w=w_{i_1}w_{i_2}\cdots w_{i_r}$ as a word in the fundamental reflections.
The weight $w(\lambda+\rho-pm\alpha)-\rho$ is given by its tuple of coefficients
with respect to the fundamental weights.
It is visually obvious, from the fact that all entries in the fourth column tuples
are nonnegative, that $m\alpha$ is a contributor. The sign of the contribution 
is given by the length of $w$.

We have discussed the immediate verifiability of the tables except
for checking that for each $w$ entry the weight $w(\lambda+\rho-pm\alpha)-\rho$
is as given.  This can be carried out by longer but routine computations (which can easily be automated).

In order to prove the theorem for a particular weight $r\omega_1$  and in characteristic
$p$,  we first find all the relevant root multiples and contributions for this weight, which may depend on $p$. Then we repeat the procedure for all contributing weights in a second iteration of the sum formula. In principle, further iterations of this process
might be expected, but for the weights being considered here 
it turns out not to be necessary; two iterations provide enough information to deduce Theorem~\ref{main}.

\subsection{Proof of Theorem~\ref{main} in detail}
\paragraph{(a)}
We may assume that $p\geq 5$. Now $V(0)$ is trivially simple, and  well known 
that $V(\omega_1)$ is simple for all $p$.
(This can also be can be checked from our tables.) 
For $p=7$, we have to check also that $V(2\omega_1)$,
and $V(3\omega_1)$ are simple. From Table~\ref{small7table1} 
we can see that no relevant root multiples are contributors.
Therefore the Jantzen sums are zero.
Assume then that $p\geq 11$. For $r<p-10$, there are no relevant roots for $r\omega_1$, so $V(r\omega_1)$ is simple.
For  $r=p-10$,\dots,$p-4$, Table~\ref{gentable} lists the relevant roots 
and shows that none is a contributor, so $V(r\omega_1)$ is simple.

\paragraph{(b)} Part (i) is obvious. For (ii), the starting point is the simplicity
of $V(\omega_1)$.  By the graph automorphism from the symmetry of the Dynkin diagram this implies that $V(\omega_6)$ is also simple. Also, in  Table~\ref{small5table1}
the only contributor for $2\omega_1$ is $\alpha=\alpha=\sum_{i=1}^5\alpha_i$, 
and the contribution is $-\chi(\omega_6)$. Hence, $\rad(V(2\omega_1))\cong L(\omega_6)$.
To prove (iii), we examine the rows of  Table~\ref{small7table2} corresponding to
$4\omega_1$, We see that the only two contributions are $\chi(\omega_1+\omega_4)$
(from $\alpha=\sum_{i=1}^6\alpha_i$) and $-\chi(2\omega_1+\omega_6)$ 
(from $\alpha=\sum_{i=1}^5\alpha_i$). We then consider the Jantzen sums
for the highest weights of these two contributions, 
by looking at Table~\ref{small7table4}. There, 
we see that $V(\omega_1+\omega_4)$ is simple.
Also, the only contribution to the Jantzen sum for $2\omega_1+\omega_6$ is
$-\chi(\omega_1+\omega_4)$, and this means that $\rad(V(2\omega_1+\omega_6))\cong V(\omega_1+\omega_4)$, and the proof of (iii) is complete.
(iv) When $p\geq 11$, the relevant roots and contributions for $(p-3)\omega_1$
are given by  Table~\ref{pm3table1}. The contributing  weights   
are $(p-9)\omega_1$, $(p-6)\omega_1+\omega_4$, $(p-8)\omega_1+\omega_2$, 
$(p-8)\omega_1+\omega_6$, and $(p-5)\omega_1+\omega_6$. The data for the 
Jantzen sums for these highest weights is given in Table~\ref{pm3table2}.
From Table~\ref{pm3table2} we see first that  $V((p-9)\omega_1)$ is simple. 
Then we see that $J((p-8)\omega_1+\omega_6)=\chi((p-9)\omega_1)$, 
which implies that $\rad(V((p-8)\omega_1+\omega_6))\cong L((p-9)\omega_1)$. 
Next, we have
\begin{equation} 
\begin{aligned}
J((p-8)\omega_1+\omega_2)&=\chi((p-8)\omega_1+\omega_6)-\chi((p-9)\omega_1)\\
&=\chi((p-8)\omega_1+\omega_6)-\Ch(\rad(V((p-8)\omega_1+\omega_6))).
\end{aligned}
\end{equation}
Hence $\rad(V((p-8)\omega_1+\omega_2))$ is a simple module isomorphic to $L((p-8)\omega_1+\omega_6)$. 
Next, we have 
\begin{equation} 
\begin{aligned}
J((p-6)\omega_1+\omega_4)&=\chi((p-8)\omega_1+\omega_2)-\chi((p-8)\omega_1+\omega_6)+\chi((p-9)\omega_1)\\
&=\chi((p-8)\omega_1+\omega_2)-J((p-8)\omega_1+\omega_2)\\
&=\chi((p-8)\omega_1+\omega_2)-\Ch(\rad(V((p-8)\omega_1+\omega_2))),
\end{aligned}
\end{equation}
which implies that $\rad(V((p-6)\omega_1+\omega_4))\cong L((p-8)\omega_1+\omega_2)$.
Continuing in the same way, we see that 
$\rad(V((p-5)\omega_1+\omega_6))\cong L((p-6)\omega_1+\omega_4)$
and that $\rad(V((p-3)\omega_1))\cong L((p-5)\omega_1+\omega_6)$. This completes the
proof of (iv) and of (b).

\paragraph{(c)} (i) The result is clear for $p=2$ and $p=3$.
(ii)For $p=5$, it is immediate from the $3\omega_1$ rows of Table~\ref{small5table2}
that the only contribution is $-\chi(0)$, so $\rad(V(\omega_3))$ is a one-dimensional
trivial module.
(iii)When $p=7$, the contributions for $5\omega_1$ can be found from 
Table~\ref{small7table2}. The contributions are $\chi(\omega_4+\omega_6)$ and $-\chi(\omega_1+2\omega_6)$. The data for a second iteration of the sum formula, applied to the
contributing weights are given in Table~\ref{small7table5}. We see that $V(\omega_4+\omega_6)$ is simple and that $\rad(V(\omega_1+2\omega_6))\cong V(\omega_4+\omega_6)$. Thus, (iii) holds.

(iv)Table~\ref{pm2table1} gives the contributions for $(p-2)\omega_1$. 
They are $-\chi((p-10)\omega_1+\omega_2)$, $\chi((p-9)\omega_1+\omega_5)$,
$-\chi((p-8)\omega_1+\omega_3)$, $\chi((p-7)\omega_1+\omega_4+\omega_6)$
and $-\chi((p-6)\omega_1+2\omega_6)$.
For each of the contributing weights, the data for a second
iteration of the sum formula are given in Table~\ref{pm2table2}.
From Table~\ref{pm2table2} we see first that $V((p-10)\omega_1+\omega_2)$ is simple. 
Then we see that $J((p-9)\omega_1+\omega_5)=\chi((p-10)\omega_1+\omega_2)$, 
which implies that $\rad(V((p-9)\omega_1+\omega_5))\cong V((p-10)\omega_1+\omega_2)$. 
Next, we have
\begin{equation} 
\begin{aligned}
J((p-8)\omega_1+\omega_3)&=\chi((p-9)\omega_1+\omega_5)-\chi((p-10)\omega_1+\omega_2)\\
&=\chi((p-9)\omega_1+\omega_5)-\Ch(\rad(V((p-9)\omega_1+\omega_5))).
\end{aligned}
\end{equation}
Hence $\rad(V((p-8)\omega_1+\omega_3))$ is a simple module isomorphic to $L((p-9)\omega_1+\omega_5)$. 
Next, we have 
\begin{equation} 
\begin{aligned}
J((p-7)\omega_1+\omega_4+\omega_6)&=\chi((p-8)\omega_1+\omega_3)-\chi((p-9)\omega_1+\omega_5)+\chi((p-10)\omega_1+\omega_2)\\
&=\chi((p-8)\omega_1+\omega_3)-J((p-8)\omega_1+\omega_3)\\
&=\chi((p-8)\omega_1+\omega_3)-\Ch(\rad(V((p-8)\omega_1+\omega_3))),
\end{aligned}
\end{equation}
which implies that $\rad(V((p-7)\omega_1+\omega_4+\omega_6))\cong L((p-8)\omega_1+\omega_3)$.
Continuing in the same way, we see that 
$\rad(V((p-6)\omega_1+2\omega_6))\cong L((p-7)\omega_1+\omega_4+\omega_6)$
and that $\rad(V((p-2)\omega_1))\cong L((p-6)\omega_1+2\omega_6)$, and the 
proof of (iv), and of (c) is complete.

\paragraph{(d)} The procedure for verifying (d) follows the  same pattern as in 
(b) and (c), and was sketched in \cite{opp}. The  steps are
as follows. Tables~\ref{small2table1}, \ref{small2table2}, \ref{small3table2}
and \ref{small5table2} yield (i). For part (ii) Table~\ref{small7table3}  shows that
$J(6\omega_1)=\chi(3\omega_6)$. The simplicity of $V(3\omega_6)$ can be checked
by a further sum formula calculation or by observing that
$V(3\omega_1)$ is simple by Table~\ref{small7table1} and the two
modules are conjugate under the graph automorphism of $G$. 
Finally, to prove (iii) we find the contributions for $V((p-1)\omega_1)$
in Table~\ref{pm1table1} and consider a second iteration the sum formula
for each of the Weyl modules for contributing weights. 
The necessary data appear in Table~\ref{pm1table2}.
In a manner completely analogous to the proof of (c)(iv), we can make the
the following sequence of inferences from Table~\ref{pm1table2}.
\begin{enumerate}
\item $V((p-11)\omega_1+2\omega_2)$ is simple.
\item $\rad(V((p-10)\omega_1+\omega_2+\omega_5))\cong L((p-11)\omega_1+2\omega_2)$.
\item $\rad(V((p-9)\omega_1+\omega_3+\omega_6))\cong L((p-10)\omega_1+\omega_2+\omega_5)$.
\item $\rad(V((p-8)\omega_1+\omega_4+2\omega_6))\cong L((p-9)\omega_1+\omega_3+\omega_6)$.
\item $\rad(V((p-7)\omega_1+3\omega_6))\cong L((p-8)\omega_1+\omega_4+2\omega_6)$.
\item $\rad(V((p-1)\omega_1))\cong L((p-7)\omega_1+3\omega_6)$.
\end{enumerate}

\paragraph{(e)} It is clear from the above discussion that (e) holds.

\begin{table}[f]
\centering
\begin{tabular}{|c|c|c|c|c|} 
\hline
$\lambda$ & $m\alpha$ & &\\
\hline
& & $\lambda+\rho-pm\alpha$ & $\beta$ \\
\cline{2-4}
$(p - 10)\omega_1$
&$[1, 2, 3, 2, 2, 1]$ & $[p - 9, 1, 1, -p + 1, 1, 1]$ & $[1, 2, 3, 1, 2, 1]$\\
\hline
& & $\lambda+\rho-pm\alpha$ & $\beta$ \\
\cline{2-4}
$(p - 9)\omega_1$
&$[1, 2, 3, 2, 2, 1]$ & $[p - 8, 1, 1, -p + 1, 1, 1]$ & $[1, 2, 2, 1, 2, 1]$\\
&$[1, 2, 3, 1, 2, 1]$ & $[p - 8, 1, -p + 1, p + 1, 1, 1]$ & $[1, 2, 2, 1, 2, 1]$\\
\hline
& & $\lambda+\rho-pm\alpha$ & $\beta$ \\
\cline{2-4}
$(p - 8)\omega_1$
&$[1, 2, 3, 2, 2, 1]$ & $[p - 7, 1, 1, -p + 1, 1, 1]$ & $[1, 2, 2, 1, 1, 1]$\\
&$[1, 2, 2, 1, 2, 1]$ & $[p - 7, -p + 1, p + 1, 1, -p + 1, 1]$ & $[1, 2, 2, 1, 1, 1]$\\
&$[1, 2, 3, 1, 2, 1]$ & $[p - 7, 1, -p + 1, p + 1, 1, 1]$ & $[1, 2, 2, 1, 1, 1]$\\
\hline
& & $\lambda+\rho-pm\alpha$ & $\beta$ \\
\cline{2-4}
$(p - 7)\omega_1$
&$[1, 2, 3, 2, 2, 1]$ & $[p - 6, 1, 1, -p + 1, 1, 1]$ & $[1, 1, 2, 1, 1, 1]$\\
&$[1, 2, 2, 1, 1, 1]$ & $[p - 6, -p + 1, 1, 1, p + 1, -p + 1]$ & $[1, 1, 2, 1, 1, 1]$\\
&$[1, 1, 2, 1, 2, 1]$ & $[-6, p + 1, 1, 1, -p + 1, 1]$ & $[1, 1, 2, 1, 1, 1]$\\
&$[1, 2, 2, 1, 2, 1]$ & $[p - 6, -p + 1, p + 1, 1, -p + 1, 1]$ & $[1, 2, 2, 1, 1, 0]$\\
&$[1, 2, 3, 1, 2, 1]$ & $[p - 6, 1, -p + 1, p + 1, 1, 1]$ & $[1, 1, 2, 1, 1, 1]$\\
\hline
& & $\lambda+\rho-pm\alpha$ & $\beta$ \\
\cline{2-4}
$(p - 6)\omega_1$
&$[1, 2, 3, 2, 2, 1]$ & $[p - 5, 1, 1, -p + 1, 1, 1]$ & $[1, 1, 1, 1, 1, 1]$\\
&$[1, 1, 2, 1, 1, 1]$ & $[-5, p + 1, -p + 1, 1, p + 1, -p + 1]$ & $[1, 1, 1, 1, 1, 1]$\\
&$[1, 2, 2, 1, 1, 1]$ & $[p - 5, -p + 1, 1, 1, p + 1, -p + 1]$ & $[1, 1, 1, 1, 1, 1]$\\
&$[1, 1, 2, 1, 2, 1]$ & $[-5, p + 1, 1, 1, -p + 1, 1]$ & $[1, 1, 1, 1, 1, 1]$\\
&$[1, 2, 2, 1, 2, 1]$ & $[p - 5, -p + 1, p + 1, 1, -p + 1, 1]$ & $[1, 1, 1, 1, 1, 1]$\\
&$[1, 2, 3, 1, 2, 1]$ & $[p - 5, 1, -p + 1, p + 1, 1, 1]$ & $[1, 1, 2, 1, 1, 0]$\\
&$[1, 2, 2, 1, 1, 0]$ & $[p - 5, -p + 1, 1, 1, 1, p + 1]$ & $[1, 1, 2, 1, 1, 0]$\\
\hline
& & $\lambda+\rho-pm\alpha$ & $\beta$ \\
\cline{2-4}
$(p - 5)\omega_1$
&$[1, 2, 3, 2, 2, 1]$ & $[p - 4, 1, 1, -p + 1, 1, 1]$ & $[1, 1, 1, 1, 1, 0]$\\
&$[1, 1, 1, 1, 1, 1]$ & $[-4, 1, p + 1, -p + 1, 1, -p + 1]$ & $[1, 1, 1, 0, 1, 1]$\\
&$[1, 1, 2, 1, 1, 1]$ & $[-4, p + 1, -p + 1, 1, p + 1, -p + 1]$ & $[1, 1, 1, 0, 1, 1]$\\
&$[1, 2, 2, 1, 1, 1]$ & $[p - 4, -p + 1, 1, 1, p + 1, -p + 1]$ & $[1, 1, 1, 0, 1, 1]$\\
&$[1, 1, 2, 1, 2, 1]$ & $[-4, p + 1, 1, 1, -p + 1, 1]$ & $[1, 1, 1, 0, 1, 1]$\\
&$[1, 2, 2, 1, 2, 1]$ & $[p - 4, -p + 1, p + 1, 1, -p + 1, 1]$ & $[1, 1, 1, 0, 1, 1]$\\
&$[1, 2, 3, 1, 2, 1]$ & $[p - 4, 1, -p + 1, p + 1, 1, 1]$ & $[1, 1, 1, 0, 1, 1]$\\
&$[1, 2, 2, 1, 1, 0]$ & $[p - 4, -p + 1, 1, 1, 1, p + 1]$ & $[1, 1, 1, 1, 1, 0]$\\
&$[1, 1, 2, 1, 1, 0]$ & $[-4, p + 1, -p + 1, 1, 1, p + 1]$ & $[1, 1, 1, 1, 1, 0]$\\
\hline
& & $\lambda+\rho-pm\alpha$ & $\beta$ \\
\cline{2-4}
$(p - 4)\omega_1$
&$[1, 2, 3, 2, 2, 1]$ & $[p - 3, 1, 1, -p + 1, 1, 1]$ & $[1, 1, 1, 1, 0, 0]$\\
&$[1, 1, 1, 1, 1, 1]$ & $[-3, 1, p + 1, -p + 1, 1, -p + 1]$ & $[1, 1, 1, 1, 0, 0]$\\
&$[1, 1, 2, 1, 1, 1]$ & $[-3, p + 1, -p + 1, 1, p + 1, -p + 1]$ & $[1, 1, 1, 1, 0, 0]$\\
&$[1, 2, 2, 1, 1, 1]$ & $[p - 3, -p + 1, 1, 1, p + 1, -p + 1]$ & $[1, 1, 1, 1, 0, 0]$\\
&$[1, 1, 2, 1, 2, 1]$ & $[-3, p + 1, 1, 1, -p + 1, 1]$ & $[1, 1, 1, 0, 1, 0]$\\
&$[1, 2, 2, 1, 2, 1]$ & $[p - 3, -p + 1, p + 1, 1, -p + 1, 1]$ & $[1, 1, 1, 0, 1, 0]$\\
&$[1, 2, 3, 1, 2, 1]$ & $[p - 3, 1, -p + 1, p + 1, 1, 1]$ & $[1, 1, 1, 0, 1, 0]$\\
&$[1, 1, 1, 0, 1, 1]$ & $[-3, 1, 1, p + 1, 1, -p + 1]$ & $[1, 1, 1, 0, 1, 0]$\\
&$[1, 2, 2, 1, 1, 0]$ & $[p - 3, -p + 1, 1, 1, 1, p + 1]$ & $[1, 1, 1, 0, 1, 0]$\\
&$[1, 1, 2, 1, 1, 0]$ & $[-3, p + 1, -p + 1, 1, 1, p + 1]$ & $[1, 1, 1, 0, 1, 0]$\\
&$[1, 1, 1, 1, 1, 0]$ & $[-3, 1, p + 1, -p + 1, -p + 1, p + 1]$ & $[1, 1, 1, 0, 1, 0]$\\
\hline
\end{tabular}
\caption{$r\omega_1$, $p-10\leq r\leq p-4$, $p\geq 11$}\label{gentable}
\end{table}
\begin{table}[f]
\centering
\begin{tabular}{|c|c|c|c|c|} 
\hline
$\lambda$ & $m\alpha$ & &\\
\hline
& & $\lambda+\rho-pm\alpha$ & $\beta$ \\
\cline{2-4}
$(p - 3)\omega_1$
&$[1, 1, 2, 1, 1, 1]$ & $[-2, p + 1, -p + 1, 1, p + 1, -p + 1]$ & $[1, 1, 1, 0, 0, 0]$\\
&$[1, 2, 2, 1, 1, 1]$ & $[p - 2, -p + 1, 1, 1, p + 1, -p + 1]$ & $[1, 1, 1, 0, 0, 0]$\\
&$[1, 2, 3, 1, 2, 1]$ & $[p - 2, 1, -p + 1, p + 1, 1, 1]$ & $[1, 1, 1, 0, 0, 0]$\\
&$[1, 1, 1, 0, 1, 1]$ & $[-2, 1, 1, p + 1, 1, -p + 1]$ & $[1, 1, 1, 0, 0, 0]$\\
&$[1, 2, 2, 1, 1, 0]$ & $[p - 2, -p + 1, 1, 1, 1, p + 1]$ & $[1, 1, 1, 0, 0, 0]$\\
&$[1, 1, 2, 1, 1, 0]$ & $[-2, p + 1, -p + 1, 1, 1, p + 1]$ & $[1, 1, 1, 0, 0, 0]$\\
&$[1, 1, 1, 0, 1, 0]$ & $[-2, 1, 1, p + 1, -p + 1, p + 1]$ & $[1, 1, 1, 0, 0, 0]$\\
&$[1, 1, 1, 1, 0, 0]$ & $[-2, 1, 1, -p + 1, p + 1, 1]$ & $[1, 1, 1, 0, 0, 0]$\\
\cline{2-4}
& & $w$ & $w(\lambda+\rho-pm\alpha)-\rho$ \\
\cline{2-4}
&$[1, 2, 3, 2, 2, 1]$ & $[1, 2, 3, 5, 4, 3, 2, 6, 5, 3, 4]$ & $[p - 9, 0, 0, 0, 0, 0]$\\
&$[1, 1, 1, 1, 1, 1]$ & $[1, 2, 3, 5, 6, 4, 2, 1]$ & $[p - 6, 0, 0, 1, 0, 0]$\\
&$[1, 1, 2, 1, 2, 1]$ & $[1, 2, 3, 5, 4, 3, 6, 5, 1]$ & $[p - 8, 1, 0, 0, 0, 0]$\\
&$[1, 2, 2, 1, 2, 1]$ & $[1, 2, 3, 5, 4, 3, 1, 6, 5, 2]$ & $[p - 8, 0, 0, 0, 0, 1]$\\
&$[1, 1, 1, 1, 1, 0]$ & $[1, 2, 3, 5, 4, 2, 1]$ & $[p - 5, 0, 0, 0, 0, 1]$\\
\hline
\end{tabular}
\caption{$(p-3)\omega_1$, $p\geq 11$}\label{pm3table1}
\end{table}
\begin{table}[f]
\centering
\begin{tabular}{|c|c|c|c|c|} 
\hline
$\lambda$ & $m\alpha$ & &\\
\hline
& & $\lambda+\rho-pm\alpha$ & $\beta$ \\
\cline{2-4}
$(p - 9)\omega_1$
&$[1, 2, 3, 2, 2, 1]$ & $[p - 8, 1, 1, -p + 1, 1, 1]$ & $[1, 2, 2, 1, 2, 1]$\\
&$[1, 2, 3, 1, 2, 1]$ & $[p - 8, 1, -p + 1, p + 1, 1, 1]$ & $[1, 2, 2, 1, 2, 1]$\\
\hline
& & $\lambda+\rho-pm\alpha$ & $\beta$ \\
\cline{2-4}
$(p - 8)\omega_1+\omega_6$
&$[1, 2, 3, 2, 2, 1]$ & $[p - 7, 1, 1, -p + 1, 1, 2]$ & $[1, 1, 2, 1, 1, 1]$\\
&$[1, 2, 2, 1, 1, 1]$ & $[p - 7, -p + 1, 1, 1, p + 1, -p + 2]$ & $[1, 1, 2, 1, 1, 1]$\\
&$[1, 1, 2, 1, 2, 1]$ & $[-7, p + 1, 1, 1, -p + 1, 2]$ & $[1, 1, 2, 1, 1, 1]$\\
&$[1, 2, 3, 1, 2, 1]$ & $[p - 7, 1, -p + 1, p + 1, 1, 2]$ & $[1, 1, 2, 1, 1, 1]$\\
\cline{2-4}
& & $w$ & $w(\lambda+\rho-pm\alpha)-\rho$ \\
\cline{2-4}
&$[1, 2, 2, 1, 2, 1]$ & $[1, 3, 2, 6, 5, 4, 3, 2, 5, 4, 3, 1, 6, 5, 2]$ & $[p - 9, 0, 0, 0, 0, 0]$\\
\hline
& & $\lambda+\rho-pm\alpha$ & $\beta$ \\
\cline{2-4}
$(p - 8)\omega_1+\omega_2$
&$[1, 2, 3, 2, 2, 1]$ & $[p - 7, 2, 1, -p + 1, 1, 1]$ & $[1, 1, 2, 1, 1, 1]$\\
&$[1, 2, 2, 1, 1, 1]$ & $[p - 7, -p + 2, 1, 1, p + 1, -p + 1]$ & $[1, 1, 2, 1, 1, 1]$\\
&$[1, 1, 2, 1, 2, 1]$ & $[-7, p + 2, 1, 1, -p + 1, 1]$ & $[1, 1, 2, 1, 1, 1]$\\
&$[1, 2, 3, 1, 2, 1]$ & $[p - 7, 2, -p + 1, p + 1, 1, 1]$ & $[1, 1, 2, 1, 1, 1]$\\
\cline{2-4}
& & $w$ & $w(\lambda+\rho-pm\alpha)-\rho$ \\
\cline{2-4}
&$[1, 2, 2, 1, 2, 1]$ & $[1, 3, 2, 5, 4, 3, 2, 5, 4, 3, 1, 6, 5, 2]$ & $[p - 9, 0, 0, 0, 0, 0]$\\
&$[1, 2, 2, 1, 1, 0]$ & $[2, 1, 3, 2, 5, 4, 3, 2, 1, 5, 4, 3, 2]$ & $[p - 8, 0, 0, 0, 0, 1]$\\
\hline
& & $\lambda+\rho-pm\alpha$ & $\beta$ \\
\cline{2-4}
$(p - 6)\omega_1+\omega_4$
&$[1, 2, 3, 2, 2, 1]$ & $[p - 5, 1, 1, -p + 2, 1, 1]$ & $[1, 1, 1, 1, 1, 0]$\\
&$[1, 1, 1, 1, 1, 1]$ & $[-5, 1, p + 1, -p + 2, 1, -p + 1]$ & $[1, 1, 1, 1, 1, 0]$\\
&$[1, 1, 2, 1, 2, 1]$ & $[-5, p + 1, 1, 2, -p + 1, 1]$ & $[1, 1, 1, 1, 1, 0]$\\
&$[1, 2, 2, 1, 2, 1]$ & $[p - 5, -p + 1, p + 1, 2, -p + 1, 1]$ & $[1, 1, 1, 1, 1, 0]$\\
&$[1, 2, 2, 1, 1, 0]$ & $[p - 5, -p + 1, 1, 2, 1, p + 1]$ & $[1, 1, 1, 1, 1, 0]$\\
&$[1, 1, 2, 1, 1, 0]$ & $[-5, p + 1, -p + 1, 2, 1, p + 1]$ & $[1, 1, 1, 1, 1, 0]$\\
\cline{2-4}
& & $w$ & $w(\lambda+\rho-pm\alpha)-\rho$ \\
\cline{2-4}
&$[1, 1, 2, 1, 1, 1]$ & $[1, 2, 4, 3, 5, 3, 2, 6, 4, 3, 1]$ & $[p - 8, 1, 0, 0, 0, 0]$\\
&$[1, 2, 2, 1, 1, 1]$ & $[1, 2, 4, 3, 5, 3, 2, 1, 6, 4, 3, 2]$ & $[p - 8, 0, 0, 0, 0, 1]$\\
&$[1, 2, 3, 1, 2, 1]$ & $[1, 2, 4, 3, 2, 1, 5, 4, 3, 2, 6, 5, 3]$ & $[p - 9, 0, 0, 0, 0, 0]$\\
\hline
& & $\lambda+\rho-pm\alpha$ & $\beta$ \\
\cline{2-4}
$(p - 5)\omega_1+\omega_6$
&$[1, 2, 3, 2, 2, 1]$ & $[p - 4, 1, 1, -p + 1, 1, 2]$ & $[1, 1, 1, 1, 1, 0]$\\
&$[1, 1, 1, 1, 1, 1]$ & $[-4, 1, p + 1, -p + 1, 1, -p + 2]$ & $[1, 1, 1, 1, 1, 0]$\\
&$[1, 1, 2, 1, 2, 1]$ & $[-4, p + 1, 1, 1, -p + 1, 2]$ & $[1, 1, 1, 1, 1, 0]$\\
&$[1, 2, 2, 1, 2, 1]$ & $[p - 4, -p + 1, p + 1, 1, -p + 1, 2]$ & $[1, 1, 1, 1, 1, 0]$\\
&$[1, 2, 2, 1, 1, 0]$ & $[p - 4, -p + 1, 1, 1, 1, p + 2]$ & $[1, 1, 1, 1, 1, 0]$\\
&$[1, 1, 2, 1, 1, 0]$ & $[-4, p + 1, -p + 1, 1, 1, p + 2]$ & $[1, 1, 1, 1, 1, 0]$\\
\cline{2-4}
& & $w$ & $w(\lambda+\rho-pm\alpha)-\rho$ \\
\cline{2-4}
&$[1, 1, 2, 1, 1, 1]$ & $[1, 2, 3, 5, 3, 2, 6, 4, 3, 1]$ & $[p - 8, 1, 0, 0, 0, 0]$\\
&$[1, 2, 2, 1, 1, 1]$ & $[1, 2, 3, 5, 3, 2, 1, 6, 4, 3, 2]$ & $[p - 8, 0, 0, 0, 0, 1]$\\
&$[1, 2, 3, 1, 2, 1]$ & $[1, 2, 3, 2, 1, 5, 4, 3, 2, 6, 5, 3]$ & $[p - 9, 0, 0, 0, 0, 0]$\\
&$[1, 1, 1, 0, 1, 1]$ & $[1, 2, 3, 5, 6, 5, 3, 2, 1]$ & $[p - 6, 0, 0, 1, 0, 0]$\\
\hline
\end{tabular}
\caption{$(p-3)\omega_1$, $p\geq 11$, second iteration.}\label{pm3table2}
\end{table}
\begin{table}[f]
\centering
\begin{tabular}{|c|c|c|c|c|} 
\hline
$\lambda$ & $m\alpha$ & &\\
\hline
& & $\lambda+\rho-pm\alpha$ & $\beta$ \\
\cline{2-4}
$(p - 2)\omega_1$
&$[1, 1, 1, 1, 1, 1]$ & $[-1, 1, p + 1, -p + 1, 1, -p + 1]$ & $[1, 1, 0, 0, 0, 0]$\\
&$[1, 2, 2, 1, 1, 1]$ & $[p - 1, -p + 1, 1, 1, p + 1, -p + 1]$ & $[1, 1, 0, 0, 0, 0]$\\
&$[1, 2, 2, 1, 2, 1]$ & $[p - 1, -p + 1, p + 1, 1, -p + 1, 1]$ & $[1, 1, 0, 0, 0, 0]$\\
&$[1, 1, 1, 0, 1, 1]$ & $[-1, 1, 1, p + 1, 1, -p + 1]$ & $[1, 1, 0, 0, 0, 0]$\\
&$[1, 2, 2, 1, 1, 0]$ & $[p - 1, -p + 1, 1, 1, 1, p + 1]$ & $[1, 1, 0, 0, 0, 0]$\\
&$[1, 1, 1, 0, 0, 0]$ & $[-1, 1, -p + 1, p + 1, p + 1, 1]$ & $[1, 1, 0, 0, 0, 0]$\\
&$[1, 1, 1, 1, 1, 0]$ & $[-1, 1, p + 1, -p + 1, -p + 1, p + 1]$ & $[1, 1, 0, 0, 0, 0]$\\
&$[1, 1, 1, 0, 1, 0]$ & $[-1, 1, 1, p + 1, -p + 1, p + 1]$ & $[1, 1, 0, 0, 0, 0]$\\
&$[1, 1, 1, 1, 0, 0]$ & $[-1, 1, 1, -p + 1, p + 1, 1]$ & $[1, 1, 0, 0, 0, 0]$\\
\cline{2-4}
& & $w$ & $w(\lambda+\rho-pm\alpha)-\rho$ \\
\cline{2-4}
&$[1, 2, 3, 2, 2, 1]$ & $[1, 2, 3, 5, 4, 3, 2, 6, 5, 3, 4]$ & $[p - 10, 1, 0, 0, 0, 0]$\\
&$[1, 1, 2, 1, 1, 1]$ & $[1, 2, 3, 5, 6, 4, 3, 1]$ & $[p - 7, 0, 0, 1, 0, 1]$\\
&$[1, 1, 2, 1, 2, 1]$ & $[1, 2, 3, 5, 4, 3, 6, 5, 1]$ & $[p - 8, 0, 1, 0, 0, 0]$\\
&$[1, 2, 3, 1, 2, 1]$ & $[1, 2, 3, 5, 4, 3, 2, 6, 5, 3]$ & $[p - 9, 0, 0, 0, 1, 0]$\\
&$[1, 1, 2, 1, 1, 0]$ & $[1, 2, 3, 5, 4, 3, 1]$ & $[p - 6, 0, 0, 0, 0, 2]$\\
\hline
\end{tabular}
\caption{$(p-2)\omega_1$, $p\geq 11$}\label{pm2table1}
\end{table}
\begin{table}[f]
\centering
\begin{tabular}{|c|c|c|c|c|} 
\hline
$\lambda$ & $m\alpha$ & &\\
\hline
& & $\lambda+\rho-pm\alpha$ & $\beta$ \\
\cline{2-4}
$(p - 10)\omega_1+\omega_2$
&$[1, 2, 3, 2, 2, 1]$ & $[p - 9, 2, 1, -p + 1, 1, 1]$ & $[1, 2, 2, 1, 1, 1]$\\
&$[1, 2, 2, 1, 2, 1]$ & $[p - 9, -p + 2, p + 1, 1, -p + 1, 1]$ & $[1, 2, 2, 1, 1, 1]$\\
&$[1, 2, 3, 1, 2, 1]$ & $[p - 9, 2, -p + 1, p + 1, 1, 1]$ & $[1, 2, 2, 1, 1, 1]$\\
\hline
& & $\lambda+\rho-pm\alpha$ & $\beta$ \\
\cline{2-4}
$(p - 9)\omega_1+\omega_5$
&$[1, 2, 3, 2, 2, 1]$ & $[p - 8, 1, 1, -p + 1, 2, 1]$ & $[1, 2, 2, 1, 1, 1]$\\
&$[1, 2, 2, 1, 2, 1]$ & $[p - 8, -p + 1, p + 1, 1, -p + 2, 1]$ & $[1, 2, 2, 1, 1, 1]$\\
&$[1, 2, 3, 1, 2, 1]$ & $[p - 8, 1, -p + 1, p + 1, 2, 1]$ & $[1, 2, 2, 1, 1, 1]$\\
\cline{2-4}
& & $w$ & $w(\lambda+\rho-pm\alpha)-\rho$ \\
\cline{2-4}
&$[1, 1, 2, 1, 2, 1]$ & $[1, 5, 3, 2, 6, 5, 4, 3, 2, 5, 4, 3, 6, 5, 1]$ & $[p - 10, 1, 0, 0, 0, 0]$\\
\hline
& & $\lambda+\rho-pm\alpha$ & $\beta$ \\
\cline{2-4}
$(p - 8)\omega_1+\omega_3$
&$[1, 2, 3, 2, 2, 1]$ & $[p - 7, 1, 2, -p + 1, 1, 1]$ & $[1, 1, 2, 1, 1, 0]$\\
&$[1, 1, 2, 1, 1, 1]$ & $[-7, p + 1, -p + 2, 1, p + 1, -p + 1]$ & $[1, 1, 2, 1, 1, 0]$\\
&$[1, 1, 2, 1, 2, 1]$ & $[-7, p + 1, 2, 1, -p + 1, 1]$ & $[1, 1, 2, 1, 1, 0]$\\
&$[1, 2, 3, 1, 2, 1]$ & $[p - 7, 1, -p + 2, p + 1, 1, 1]$ & $[1, 1, 2, 1, 1, 0]$\\
&$[1, 2, 2, 1, 1, 0]$ & $[p - 7, -p + 1, 2, 1, 1, p + 1]$ & $[1, 1, 2, 1, 1, 0]$\\
\cline{2-4}
& & $w$ & $w(\lambda+\rho-pm\alpha)-\rho$ \\
\cline{2-4}
&$[1, 2, 2, 1, 1, 1]$ & $[1, 3, 2, 4, 3, 5, 3, 2, 1, 6, 4, 3, 2]$ & $[p - 9, 0, 0, 0, 1, 0]$\\
&$[1, 2, 2, 1, 2, 1]$ & $[1, 3, 2, 5, 4, 3, 2, 5, 4, 3, 1, 6, 5, 2]$ & $[p - 10, 1, 0, 0, 0, 0]$\\
\hline
& & $\lambda+\rho-pm\alpha$ & $\beta$ \\
\cline{2-4}
$(p - 7)\omega_1+\omega_4+\omega_6$
&$[1, 2, 3, 2, 2, 1]$ & $[p - 6, 1, 1, -p + 2, 1, 2]$ & $[1, 1, 2, 1, 1, 0]$\\
&$[1, 1, 2, 1, 1, 1]$ & $[-6, p + 1, -p + 1, 2, p + 1, -p + 2]$ & $[1, 1, 2, 1, 1, 0]$\\
&$[1, 1, 2, 1, 2, 1]$ & $[-6, p + 1, 1, 2, -p + 1, 2]$ & $[1, 1, 2, 1, 1, 0]$\\
&$[1, 2, 3, 1, 2, 1]$ & $[p - 6, 1, -p + 1, p + 2, 1, 2]$ & $[1, 1, 2, 1, 1, 0]$\\
&$[1, 2, 2, 1, 1, 0]$ & $[p - 6, -p + 1, 1, 2, 1, p + 2]$ & $[1, 1, 2, 1, 1, 0]$\\
\cline{2-4}
& & $w$ & $w(\lambda+\rho-pm\alpha)-\rho$ \\
\cline{2-4}
&$[1, 1, 1, 1, 1, 1]$ & $[1, 2, 4, 3, 6, 5, 3, 6, 4, 2, 1]$ & $[p - 8, 0, 1, 0, 0, 0]$\\
&$[1, 2, 2, 1, 1, 1]$ & $[1, 2, 4, 3, 5, 3, 2, 1, 6, 4, 3, 2]$ & $[p - 9, 0, 0, 0, 1, 0]$\\
&$[1, 2, 2, 1, 2, 1]$ & $[1, 2, 5, 4, 3, 2, 5, 4, 3, 1, 6, 5, 2]$ & $[p - 10, 1, 0, 0, 0, 0]$\\
\hline
& & $\lambda+\rho-pm\alpha$ & $\beta$ \\
\cline{2-4}
$(p - 6)\omega_1+2\omega_6$
&$[1, 2, 3, 2, 2, 1]$ & $[p - 5, 1, 1, -p + 1, 1, 3]$ & $[1, 1, 2, 1, 1, 0]$\\
&$[1, 1, 2, 1, 1, 1]$ & $[-5, p + 1, -p + 1, 1, p + 1, -p + 3]$ & $[1, 1, 2, 1, 1, 0]$\\
&$[1, 1, 2, 1, 2, 1]$ & $[-5, p + 1, 1, 1, -p + 1, 3]$ & $[1, 1, 2, 1, 1, 0]$\\
&$[1, 2, 3, 1, 2, 1]$ & $[p - 5, 1, -p + 1, p + 1, 1, 3]$ & $[1, 1, 2, 1, 1, 0]$\\
&$[1, 2, 2, 1, 1, 0]$ & $[p - 5, -p + 1, 1, 1, 1, p + 3]$ & $[1, 1, 2, 1, 1, 0]$\\
\cline{2-4}
& & $w$ & $w(\lambda+\rho-pm\alpha)-\rho$ \\
\cline{2-4}
&$[1, 1, 1, 1, 1, 1]$ & $[1, 2, 3, 6, 5, 3, 6, 4, 2, 1]$ & $[p - 8, 0, 1, 0, 0, 0]$\\
&$[1, 2, 2, 1, 1, 1]$ & $[1, 2, 3, 5, 3, 2, 1, 6, 4, 3, 2]$ & $[p - 9, 0, 0, 0, 1, 0]$\\
&$[1, 2, 2, 1, 2, 1]$ & $[1, 2, 5, 3, 2, 5, 4, 3, 1, 6, 5, 2]$ & $[p - 10, 1, 0, 0, 0, 0]$\\
&$[1, 1, 1, 0, 1, 1]$ & $[1, 2, 3, 5, 6, 5, 3, 2, 1]$ & $[p - 7, 0, 0, 1, 0, 1]$\\
\hline
\end{tabular}
\caption{$(p-2)\omega_1$, $p\geq 11$, second iteration}\label{pm2table2}
\end{table}
\begin{table}[f]
\centering
\begin{tabular}{|c|c|c|c|c|} 
\hline
$\lambda$ & $m\alpha$ & &\\
\hline
& & $\lambda+\rho-pm\alpha$ & $\beta$ \\
\cline{2-4}
$(p - 1)\omega_1$
&$[1, 1, 1, 1, 1, 1]$ & $[0, 1, p + 1, -p + 1, 1, -p + 1]$ & $[1, 0, 0, 0, 0, 0]$\\
&$[1, 1, 2, 1, 1, 1]$ & $[0, p + 1, -p + 1, 1, p + 1, -p + 1]$ & $[1, 0, 0, 0, 0, 0]$\\
&$[1, 1, 2, 1, 2, 1]$ & $[0, p + 1, 1, 1, -p + 1, 1]$ & $[1, 0, 0, 0, 0, 0]$\\
&$[1, 1, 1, 0, 1, 1]$ & $[0, 1, 1, p + 1, 1, -p + 1]$ & $[1, 0, 0, 0, 0, 0]$\\
&$[1, 1, 0, 0, 0, 0]$ & $[0, -p + 1, p + 1, 1, 1, 1]$ & $[1, 0, 0, 0, 0, 0]$\\
&$[1, 1, 2, 1, 1, 0]$ & $[0, p + 1, -p + 1, 1, 1, p + 1]$ & $[1, 0, 0, 0, 0, 0]$\\
&$[1, 1, 1, 0, 0, 0]$ & $[0, 1, -p + 1, p + 1, p + 1, 1]$ & $[1, 0, 0, 0, 0, 0]$\\
&$[1, 1, 1, 1, 1, 0]$ & $[0, 1, p + 1, -p + 1, -p + 1, p + 1]$ & $[1, 0, 0, 0, 0, 0]$\\
&$[1, 1, 1, 0, 1, 0]$ & $[0, 1, 1, p + 1, -p + 1, p + 1]$ & $[1, 0, 0, 0, 0, 0]$\\
&$[1, 1, 1, 1, 0, 0]$ & $[0, 1, 1, -p + 1, p + 1, 1]$ & $[1, 0, 0, 0, 0, 0]$\\
\cline{2-4}
& & $w$ & $w(\lambda+\rho-pm\alpha)-\rho$ \\
\cline{2-4}
&$[1, 2, 3, 2, 2, 1]$ & $[1, 2, 3, 5, 4, 3, 2, 6, 5, 3, 4]$ & $[p - 11, 2, 0, 0, 0, 0]$\\
&$[1, 2, 2, 1, 1, 1]$ & $[1, 2, 3, 5, 6, 4, 3, 2]$ & $[p - 8, 0, 0, 1, 0, 2]$\\
&$[1, 2, 2, 1, 2, 1]$ & $[1, 2, 3, 5, 4, 3, 6, 5, 2]$ & $[p - 9, 0, 1, 0, 0, 1]$\\
&$[1, 2, 3, 1, 2, 1]$ & $[1, 2, 3, 5, 4, 3, 2, 6, 5, 3]$ & $[p - 10, 1, 0, 0, 1, 0]$\\
&$[1, 2, 2, 1, 1, 0]$ & $[1, 2, 3, 5, 4, 3, 2]$ & $[p - 7, 0, 0, 0, 0, 3]$\\
\hline
\end{tabular}
\caption{$(p-1)\omega_1$, $p\geq 11$}\label{pm1table1}
\end{table}
\begin{table}[f]
\centering
\begin{tabular}{|c|c|c|c|c|} 
\hline
$\lambda$ & $m\alpha$ & &\\
\hline
& & $\lambda+\rho-pm\alpha$ & $\beta$ \\
\cline{2-4}
$(p - 11)\omega_1+2\omega_2$
&$[1, 2, 3, 2, 2, 1]$ & $[p - 10, 3, 1, -p + 1, 1, 1]$ & $[1, 2, 2, 1, 1, 0]$\\
&$[1, 2, 2, 1, 1, 1]$ & $[p - 10, -p + 3, 1, 1, p + 1, -p + 1]$ & $[1, 2, 2, 1, 1, 0]$\\
&$[1, 2, 2, 1, 2, 1]$ & $[p - 10, -p + 3, p + 1, 1, -p + 1, 1]$ & $[1, 2, 2, 1, 1, 0]$\\
&$[1, 2, 3, 1, 2, 1]$ & $[p - 10, 3, -p + 1, p + 1, 1, 1]$ & $[1, 2, 2, 1, 1, 0]$\\
\hline
& & $\lambda+\rho-pm\alpha$ & $\beta$ \\
\cline{2-4}
$(p - 10)\omega_1+\omega_2+\omega_5$
&$[1, 2, 3, 2, 2, 1]$ & $[p - 9, 2, 1, -p + 1, 2, 1]$ & $[1, 2, 2, 1, 1, 0]$\\
&$[1, 2, 2, 1, 1, 1]$ & $[p - 9, -p + 2, 1, 1, p + 2, -p + 1]$ & $[1, 2, 2, 1, 1, 0]$\\
&$[1, 2, 2, 1, 2, 1]$ & $[p - 9, -p + 2, p + 1, 1, -p + 2, 1]$ & $[1, 2, 2, 1, 1, 0]$\\
&$[1, 2, 3, 1, 2, 1]$ & $[p - 9, 2, -p + 1, p + 1, 2, 1]$ & $[1, 2, 2, 1, 1, 0]$\\
\cline{2-4}
& & $w$ & $w(\lambda+\rho-pm\alpha)-\rho$ \\
\cline{2-4}
&$[1, 1, 2, 1, 2, 1]$ & $[1, 5, 3, 2, 6, 5, 4, 3, 2, 5, 4, 3, 6, 5, 1]$ & $[p - 11, 2, 0, 0, 0, 0]$\\
\hline
& & $\lambda+\rho-pm\alpha$ & $\beta$ \\
\cline{2-4}
$(p - 9)\omega_1+\omega_3+\omega_6$
&$[1, 2, 3, 2, 2, 1]$ & $[p - 8, 1, 2, -p + 1, 1, 2]$ & $[1, 2, 2, 1, 1, 0]$\\
&$[1, 2, 2, 1, 1, 1]$ & $[p - 8, -p + 1, 2, 1, p + 1, -p + 2]$ & $[1, 2, 2, 1, 1, 0]$\\
&$[1, 2, 2, 1, 2, 1]$ & $[p - 8, -p + 1, p + 2, 1, -p + 1, 2]$ & $[1, 2, 2, 1, 1, 0]$\\
&$[1, 2, 3, 1, 2, 1]$ & $[p - 8, 1, -p + 2, p + 1, 1, 2]$ & $[1, 2, 2, 1, 1, 0]$\\
\cline{2-4}
& & $w$ & $w(\lambda+\rho-pm\alpha)-\rho$ \\
\cline{2-4}
&$[1, 1, 2, 1, 1, 1]$ & $[1, 3, 2, 4, 3, 6, 5, 3, 2, 6, 4, 3, 1]$ & $[p - 10, 1, 0, 0, 1, 0]$\\
&$[1, 1, 2, 1, 2, 1]$ & $[1, 3, 2, 6, 5, 4, 3, 2, 5, 4, 3, 6, 5, 1]$ & $[p - 11, 2, 0, 0, 0, 0]$\\
\hline
& & $\lambda+\rho-pm\alpha$ & $\beta$ \\
\cline{2-4}
$(p - 8)\omega_1+\omega_4+2\omega_6$
&$[1, 2, 3, 2, 2, 1]$ & $[p - 7, 1, 1, -p + 2, 1, 3]$ & $[1, 2, 2, 1, 1, 0]$\\
&$[1, 2, 2, 1, 1, 1]$ & $[p - 7, -p + 1, 1, 2, p + 1, -p + 3]$ & $[1, 2, 2, 1, 1, 0]$\\
&$[1, 2, 2, 1, 2, 1]$ & $[p - 7, -p + 1, p + 1, 2, -p + 1, 3]$ & $[1, 2, 2, 1, 1, 0]$\\
&$[1, 2, 3, 1, 2, 1]$ & $[p - 7, 1, -p + 1, p + 2, 1, 3]$ & $[1, 2, 2, 1, 1, 0]$\\
\cline{2-4}
& & $w$ & $w(\lambda+\rho-pm\alpha)-\rho$ \\
\cline{2-4}
&$[1, 1, 1, 1, 1, 1]$ & $[1, 2, 4, 3, 6, 5, 3, 6, 4, 2, 1]$ & $[p - 9, 0, 1, 0, 0, 1]$\\
&$[1, 1, 2, 1, 1, 1]$ & $[1, 2, 4, 3, 6, 5, 3, 2, 6, 4, 3, 1]$ & $[p - 10, 1, 0, 0, 1, 0]$\\
&$[1, 1, 2, 1, 2, 1]$ & $[1, 2, 6, 5, 4, 3, 2, 5, 4, 3, 6, 5, 1]$ & $[p - 11, 2, 0, 0, 0, 0]$\\
\hline
& & $\lambda+\rho-pm\alpha$ & $\beta$ \\
\cline{2-4}
$(p - 7)\omega_1+3\omega_6$
&$[1, 2, 3, 2, 2, 1]$ & $[p - 6, 1, 1, -p + 1, 1, 4]$ & $[1, 2, 2, 1, 1, 0]$\\
&$[1, 2, 2, 1, 1, 1]$ & $[p - 6, -p + 1, 1, 1, p + 1, -p + 4]$ & $[1, 2, 2, 1, 1, 0]$\\
&$[1, 2, 2, 1, 2, 1]$ & $[p - 6, -p + 1, p + 1, 1, -p + 1, 4]$ & $[1, 2, 2, 1, 1, 0]$\\
&$[1, 2, 3, 1, 2, 1]$ & $[p - 6, 1, -p + 1, p + 1, 1, 4]$ & $[1, 2, 2, 1, 1, 0]$\\
\cline{2-4}
& & $w$ & $w(\lambda+\rho-pm\alpha)-\rho$ \\
\cline{2-4}
&$[1, 1, 1, 1, 1, 1]$ & $[1, 2, 3, 6, 5, 3, 6, 4, 2, 1]$ & $[p - 9, 0, 1, 0, 0, 1]$\\
&$[1, 1, 2, 1, 1, 1]$ & $[1, 2, 3, 6, 5, 3, 2, 6, 4, 3, 1]$ & $[p - 10, 1, 0, 0, 1, 0]$\\
&$[1, 1, 2, 1, 2, 1]$ & $[1, 2, 6, 5, 3, 2, 5, 4, 3, 6, 5, 1]$ & $[p - 11, 2, 0, 0, 0, 0]$\\
&$[1, 1, 1, 0, 1, 1]$ & $[1, 2, 3, 5, 6, 5, 3, 2, 1]$ & $[p - 8, 0, 0, 1, 0, 2]$\\
\hline
\end{tabular}
\caption{$(p-1)\omega_1$, $p\geq 11$, second iteration}\label{pm1table2}
\end{table}
\begin{table}[f]
\centering
\begin{tabular}{|c|c|c|c|c|} 
\hline
$\lambda$ & $m\alpha$ & &\\
\hline
& & $\lambda+\rho-pm\alpha$ & $\beta$ \\
\cline{2-4}
$\omega_1$
&$[1, 2, 3, 2, 2, 1]$ & $[2, 1, 1, -1, 1, 1]$ & $[0, 0, 1, 1, 0, 0]$\\
&$[2, 4, 6, 4, 4, 2]$ & $[2, 1, 1, -3, 1, 1]$ & $[0, 0, 1, 1, 1, 1]$\\
&$[3, 6, 9, 6, 6, 3]$ & $[2, 1, 1, -5, 1, 1]$ & $[1, 2, 3, 2, 2, 1]$\\
&$[4, 8, 12, 8, 8, 4]$ & $[2, 1, 1, -7, 1, 1]$ & $[1, 1, 2, 1, 1, 1]$\\
&$[5, 10, 15, 10, 10, 5]$ & $[2, 1, 1, -9, 1, 1]$ & $[1, 2, 2, 1, 2, 1]$\\
&$[0, 0, 1, 1, 1, 1]$ & $[2, 3, 1, -1, 1, -1]$ & $[0, 0, 1, 1, 1, 1]$\\
&$[0, 1, 1, 1, 1, 1]$ & $[4, -1, 3, -1, 1, -1]$ & $[0, 0, 0, 0, 1, 1]$\\
&$[0, 2, 2, 2, 2, 2]$ & $[6, -3, 5, -3, 1, -3]$ & $[0, 0, 1, 1, 1, 1]$\\
&$[1, 1, 1, 1, 1, 1]$ & $[0, 1, 3, -1, 1, -1]$ & $[0, 0, 0, 0, 1, 1]$\\
&$[2, 2, 2, 2, 2, 2]$ & $[-2, 1, 5, -3, 1, -3]$ & $[0, 0, 1, 1, 1, 1]$\\
&$[3, 3, 3, 3, 3, 3]$ & $[-4, 1, 7, -5, 1, -5]$ & $[1, 1, 1, 0, 1, 1]$\\
&$[0, 1, 2, 1, 1, 1]$ & $[4, 1, -1, 1, 3, -1]$ & $[0, 1, 1, 0, 0, 0]$\\
&$[0, 2, 4, 2, 2, 2]$ & $[6, 1, -3, 1, 5, -3]$ & $[0, 0, 1, 1, 1, 1]$\\
&$[1, 1, 2, 1, 1, 1]$ & $[0, 3, -1, 1, 3, -1]$ & $[1, 0, 0, 0, 0, 0]$\\
&$[2, 2, 4, 2, 2, 2]$ & $[-2, 5, -3, 1, 5, -3]$ & $[0, 0, 1, 1, 1, 1]$\\
&$[3, 3, 6, 3, 3, 3]$ & $[-4, 7, -5, 1, 7, -5]$ & $[0, 1, 2, 1, 1, 1]$\\
&$[1, 2, 2, 1, 1, 1]$ & $[2, -1, 1, 1, 3, -1]$ & $[0, 1, 1, 0, 0, 0]$\\
&$[2, 4, 4, 2, 2, 2]$ & $[2, -3, 1, 1, 5, -3]$ & $[0, 1, 1, 0, 1, 1]$\\
&$[3, 6, 6, 3, 3, 3]$ & $[2, -5, 1, 1, 7, -5]$ & $[0, 1, 2, 1, 1, 1]$\\
&$[4, 8, 8, 4, 4, 4]$ & $[2, -7, 1, 1, 9, -7]$ & $[1, 1, 2, 1, 1, 1]$\\
&$[0, 1, 2, 1, 2, 1]$ & $[4, 1, 1, 1, -1, 1]$ & $[0, 0, 0, 0, 1, 1]$\\
&$[0, 2, 4, 2, 4, 2]$ & $[6, 1, 1, 1, -3, 1]$ & $[0, 0, 1, 1, 1, 1]$\\
&$[0, 3, 6, 3, 6, 3]$ & $[8, 1, 1, 1, -5, 1]$ & $[0, 1, 2, 1, 1, 1]$\\
&$[1, 1, 2, 1, 2, 1]$ & $[0, 3, 1, 1, -1, 1]$ & $[0, 0, 0, 0, 1, 1]$\\
&$[2, 2, 4, 2, 4, 2]$ & $[-2, 5, 1, 1, -3, 1]$ & $[0, 0, 1, 1, 1, 1]$\\
&$[3, 3, 6, 3, 6, 3]$ & $[-4, 7, 1, 1, -5, 1]$ & $[1, 1, 1, 0, 1, 1]$\\
&$[4, 4, 8, 4, 8, 4]$ & $[-6, 9, 1, 1, -7, 1]$ & $[1, 1, 2, 1, 1, 1]$\\
&$[1, 2, 2, 1, 2, 1]$ & $[2, -1, 3, 1, -1, 1]$ & $[0, 0, 0, 0, 1, 1]$\\
&$[2, 4, 4, 2, 4, 2]$ & $[2, -3, 5, 1, -3, 1]$ & $[0, 1, 1, 0, 1, 1]$\\
&$[3, 6, 6, 3, 6, 3]$ & $[2, -5, 7, 1, -5, 1]$ & $[1, 1, 1, 0, 1, 1]$\\
&$[4, 8, 8, 4, 8, 4]$ & $[2, -7, 9, 1, -7, 1]$ & $[1, 2, 2, 1, 1, 0]$\\
&$[1, 2, 3, 1, 2, 1]$ & $[2, 1, -1, 3, 1, 1]$ & $[0, 1, 1, 0, 0, 0]$\\
&$[2, 4, 6, 2, 4, 2]$ & $[2, 1, -3, 5, 1, 1]$ & $[0, 1, 1, 0, 1, 1]$\\
&$[3, 6, 9, 3, 6, 3]$ & $[2, 1, -5, 7, 1, 1]$ & $[0, 1, 2, 1, 1, 1]$\\
&$[4, 8, 12, 4, 8, 4]$ & $[2, 1, -7, 9, 1, 1]$ & $[1, 1, 2, 1, 1, 1]$\\
&$[5, 10, 15, 5, 10, 5]$ & $[2, 1, -9, 11, 1, 1]$ & $[1, 2, 2, 1, 2, 1]$\\
\hline
\end{tabular}
\caption{$p=2$, relevant root multiples for $\omega_1$, Part 1.}\label{small2table1}
\end{table}
\begin{table}[f]
\centering
\begin{tabular}{|c|c|c|c|c|} 
\hline
$\lambda$ & $m\alpha$ & &\\
\hline
& & $\lambda+\rho-pm\alpha$ & $\beta$ \\
\cline{2-4}
$\omega_1$
&$[0, 0, 1, 0, 1, 1]$ & $[2, 3, -1, 3, 1, -1]$ & $[0, 0, 0, 0, 1, 1]$\\
&$[0, 1, 1, 0, 1, 1]$ & $[4, -1, 1, 3, 1, -1]$ & $[0, 0, 0, 0, 1, 1]$\\
&$[1, 1, 1, 0, 1, 1]$ & $[0, 1, 1, 3, 1, -1]$ & $[0, 0, 0, 0, 1, 1]$\\
&$[2, 2, 2, 0, 2, 2]$ & $[-2, 1, 1, 5, 1, -3]$ & $[0, 1, 1, 0, 1, 1]$\\
&$[1, 2, 2, 1, 1, 0]$ & $[2, -1, 1, 1, 1, 3]$ & $[0, 1, 1, 0, 0, 0]$\\
&$[2, 4, 4, 2, 2, 0]$ & $[2, -3, 1, 1, 1, 5]$ & $[1, 2, 2, 1, 1, 0]$\\
&$[3, 6, 6, 3, 3, 0]$ & $[2, -5, 1, 1, 1, 7]$ & $[1, 1, 1, 1, 1, 0]$\\
&$[1, 1, 0, 0, 0, 0]$ & $[0, -1, 3, 1, 1, 1]$ & $[1, 0, 0, 0, 0, 0]$\\
&$[1, 1, 2, 1, 1, 0]$ & $[0, 3, -1, 1, 1, 3]$ & $[1, 0, 0, 0, 0, 0]$\\
&$[2, 2, 4, 2, 2, 0]$ & $[-2, 5, -3, 1, 1, 5]$ & $[1, 1, 1, 0, 0, 0]$\\
&$[3, 3, 6, 3, 3, 0]$ & $[-4, 7, -5, 1, 1, 7]$ & $[1, 1, 1, 1, 1, 0]$\\
&$[1, 1, 1, 0, 0, 0]$ & $[0, 1, -1, 3, 3, 1]$ & $[1, 0, 0, 0, 0, 0]$\\
&$[1, 1, 1, 1, 1, 0]$ & $[0, 1, 3, -1, -1, 3]$ & $[1, 0, 0, 0, 0, 0]$\\
&$[2, 2, 2, 2, 2, 0]$ & $[-2, 1, 5, -3, -3, 5]$ & $[0, 1, 1, 1, 1, 0]$\\
&$[1, 1, 1, 0, 1, 0]$ & $[0, 1, 1, 3, -1, 3]$ & $[1, 0, 0, 0, 0, 0]$\\
&$[2, 2, 2, 0, 2, 0]$ & $[-2, 1, 1, 5, -3, 5]$ & $[1, 1, 1, 0, 0, 0]$\\
&$[1, 1, 1, 1, 0, 0]$ & $[0, 1, 1, -1, 3, 1]$ & $[1, 0, 0, 0, 0, 0]$\\
&$[2, 2, 2, 2, 0, 0]$ & $[-2, 1, 1, -3, 5, 1]$ & $[1, 1, 1, 0, 0, 0]$\\
&$[0, 1, 2, 1, 1, 0]$ & $[4, 1, -1, 1, 1, 3]$ & $[0, 1, 1, 0, 0, 0]$\\
&$[0, 2, 4, 2, 2, 0]$ & $[6, 1, -3, 1, 1, 5]$ & $[0, 1, 1, 1, 1, 0]$\\
&$[0, 1, 1, 1, 1, 0]$ & $[4, -1, 3, -1, -1, 3]$ & $[0, 1, 1, 1, 1, 0]$\\
&$[0, 1, 1, 0, 1, 0]$ & $[4, -1, 1, 3, -1, 3]$ & $[0, 1, 1, 0, 0, 0]$\\
&$[0, 1, 1, 1, 0, 0]$ & $[4, -1, 1, -1, 3, 1]$ & $[0, 1, 1, 0, 0, 0]$\\
&$[0, 0, 1, 1, 1, 0]$ & $[2, 3, 1, -1, -1, 3]$ & $[0, 0, 1, 0, 1, 0]$\\
\hline
\end{tabular}
\caption{$p=2$, relevant root multiples for $\omega_1$, Part 2.}\label{small2table2}
\end{table}
\begin{table}[f]
\centering
\begin{tabular}{|c|c|c|c|c|} 
\hline
$\lambda$ & $m\alpha$ & &\\
\hline
& & $\lambda+\rho-pm\alpha$ & $\beta$ \\
\cline{2-4}
$\omega_1$
&$[1, 2, 3, 2, 2, 1]$ & $[2, 1, 1, -2, 1, 1]$ & $[0, 1, 1, 1, 0, 0]$\\
&$[2, 4, 6, 4, 4, 2]$ & $[2, 1, 1, -5, 1, 1]$ & $[1, 2, 3, 2, 2, 1]$\\
&$[3, 6, 9, 6, 6, 3]$ & $[2, 1, 1, -8, 1, 1]$ & $[1, 2, 2, 1, 1, 1]$\\
&$[0, 0, 1, 1, 1, 1]$ & $[2, 4, 1, -2, 1, -2]$ & $[0, 0, 1, 0, 1, 1]$\\
&$[0, 1, 1, 1, 1, 1]$ & $[5, -2, 4, -2, 1, -2]$ & $[0, 1, 1, 1, 0, 0]$\\
&$[1, 1, 1, 1, 1, 1]$ & $[-1, 1, 4, -2, 1, -2]$ & $[1, 1, 0, 0, 0, 0]$\\
&$[2, 2, 2, 2, 2, 2]$ & $[-4, 1, 7, -5, 1, -5]$ & $[1, 1, 1, 0, 1, 1]$\\
&$[0, 1, 2, 1, 1, 1]$ & $[5, 1, -2, 1, 4, -2]$ & $[0, 1, 2, 1, 1, 1]$\\
&$[1, 1, 2, 1, 1, 1]$ & $[-1, 4, -2, 1, 4, -2]$ & $[0, 0, 1, 0, 1, 1]$\\
&$[2, 2, 4, 2, 2, 2]$ & $[-4, 7, -5, 1, 7, -5]$ & $[0, 1, 2, 1, 1, 1]$\\
&$[1, 2, 2, 1, 1, 1]$ & $[2, -2, 1, 1, 4, -2]$ & $[1, 1, 0, 0, 0, 0]$\\
&$[2, 4, 4, 2, 2, 2]$ & $[2, -5, 1, 1, 7, -5]$ & $[0, 1, 2, 1, 1, 1]$\\
&$[0, 1, 2, 1, 2, 1]$ & $[5, 1, 1, 1, -2, 1]$ & $[0, 0, 1, 0, 1, 1]$\\
&$[0, 2, 4, 2, 4, 2]$ & $[8, 1, 1, 1, -5, 1]$ & $[0, 1, 2, 1, 1, 1]$\\
&$[1, 1, 2, 1, 2, 1]$ & $[-1, 4, 1, 1, -2, 1]$ & $[0, 0, 1, 0, 1, 1]$\\
&$[2, 2, 4, 2, 4, 2]$ & $[-4, 7, 1, 1, -5, 1]$ & $[1, 1, 1, 0, 1, 1]$\\
&$[1, 2, 2, 1, 2, 1]$ & $[2, -2, 4, 1, -2, 1]$ & $[1, 1, 0, 0, 0, 0]$\\
&$[2, 4, 4, 2, 4, 2]$ & $[2, -5, 7, 1, -5, 1]$ & $[1, 1, 1, 0, 1, 1]$\\
&$[3, 6, 6, 3, 6, 3]$ & $[2, -8, 10, 1, -8, 1]$ & $[1, 2, 2, 1, 1, 1]$\\
&$[1, 2, 3, 1, 2, 1]$ & $[2, 1, -2, 4, 1, 1]$ & $[0, 0, 1, 0, 1, 1]$\\
&$[2, 4, 6, 2, 4, 2]$ & $[2, 1, -5, 7, 1, 1]$ & $[0, 1, 2, 1, 1, 1]$\\
&$[3, 6, 9, 3, 6, 3]$ & $[2, 1, -8, 10, 1, 1]$ & $[1, 2, 2, 1, 1, 1]$\\
&$[0, 1, 1, 0, 1, 1]$ & $[5, -2, 1, 4, 1, -2]$ & $[0, 0, 1, 0, 1, 1]$\\
&$[1, 1, 1, 0, 1, 1]$ & $[-1, 1, 1, 4, 1, -2]$ & $[0, 0, 1, 0, 1, 1]$\\
&$[1, 2, 2, 1, 1, 0]$ & $[2, -2, 1, 1, 1, 4]$ & $[1, 1, 0, 0, 0, 0]$\\
&$[2, 4, 4, 2, 2, 0]$ & $[2, -5, 1, 1, 1, 7]$ & $[1, 1, 1, 1, 1, 0]$\\
&$[1, 1, 2, 1, 1, 0]$ & $[-1, 4, -2, 1, 1, 4]$ & $[0, 0, 1, 1, 1, 0]$\\
&$[2, 2, 4, 2, 2, 0]$ & $[-4, 7, -5, 1, 1, 7]$ & $[1, 1, 1, 1, 1, 0]$\\
&$[1, 1, 1, 0, 0, 0]$ & $[-1, 1, -2, 4, 4, 1]$ & $[1, 1, 0, 0, 0, 0]$\\
&$[1, 1, 1, 1, 1, 0]$ & $[-1, 1, 4, -2, -2, 4]$ & $[1, 1, 0, 0, 0, 0]$\\
&$[1, 1, 1, 0, 1, 0]$ & $[-1, 1, 1, 4, -2, 4]$ & $[1, 1, 0, 0, 0, 0]$\\
&$[1, 1, 1, 1, 0, 0]$ & $[-1, 1, 1, -2, 4, 1]$ & $[1, 1, 0, 0, 0, 0]$\\
&$[0, 1, 2, 1, 1, 0]$ & $[5, 1, -2, 1, 1, 4]$ & $[0, 1, 1, 0, 1, 0]$\\
&$[0, 1, 1, 1, 1, 0]$ & $[5, -2, 4, -2, -2, 4]$ & $[0, 1, 1, 0, 1, 0]$\\
\hline
\end{tabular}
\caption{$p=3$, relevant root multiples for $\omega_1$}\label{small3table1}
\end{table}
\begin{table}[f]
\centering
\begin{tabular}{|c|c|c|c|c|} 
\hline
$\lambda$ & $m\alpha$ & &\\
\hline
& & $\lambda+\rho-pm\alpha$ & $\beta$ \\
\cline{2-4}
$2\omega_1$
&$[1, 2, 3, 2, 2, 1]$ & $[3, 1, 1, -2, 1, 1]$ & $[0, 1, 1, 1, 0, 0]$\\
&$[2, 4, 6, 4, 4, 2]$ & $[3, 1, 1, -5, 1, 1]$ & $[0, 1, 2, 1, 1, 1]$\\
&$[3, 6, 9, 6, 6, 3]$ & $[3, 1, 1, -8, 1, 1]$ & $[1, 1, 2, 1, 1, 1]$\\
&$[4, 8, 12, 8, 8, 4]$ & $[3, 1, 1, -11, 1, 1]$ & $[1, 2, 3, 1, 2, 1]$\\
&$[0, 0, 1, 1, 1, 1]$ & $[3, 4, 1, -2, 1, -2]$ & $[0, 0, 1, 0, 1, 1]$\\
&$[0, 1, 1, 1, 1, 1]$ & $[6, -2, 4, -2, 1, -2]$ & $[0, 1, 1, 1, 0, 0]$\\
&$[1, 1, 1, 1, 1, 1]$ & $[0, 1, 4, -2, 1, -2]$ & $[1, 0, 0, 0, 0, 0]$\\
&$[2, 2, 2, 2, 2, 2]$ & $[-3, 1, 7, -5, 1, -5]$ & $[1, 1, 1, 1, 0, 0]$\\
&$[0, 1, 2, 1, 1, 1]$ & $[6, 1, -2, 1, 4, -2]$ & $[0, 1, 2, 1, 1, 1]$\\
&$[1, 1, 2, 1, 1, 1]$ & $[0, 4, -2, 1, 4, -2]$ & $[0, 0, 1, 0, 1, 1]$\\
&$[2, 2, 4, 2, 2, 2]$ & $[-3, 7, -5, 1, 7, -5]$ & $[0, 1, 2, 1, 1, 1]$\\
&$[1, 2, 2, 1, 1, 1]$ & $[3, -2, 1, 1, 4, -2]$ & $[0, 1, 1, 1, 0, 0]$\\
&$[2, 4, 4, 2, 2, 2]$ & $[3, -5, 1, 1, 7, -5]$ & $[0, 1, 2, 1, 1, 1]$\\
&$[3, 6, 6, 3, 3, 3]$ & $[3, -8, 1, 1, 10, -8]$ & $[1, 1, 2, 1, 1, 1]$\\
&$[0, 1, 2, 1, 2, 1]$ & $[6, 1, 1, 1, -2, 1]$ & $[0, 0, 1, 0, 1, 1]$\\
&$[0, 2, 4, 2, 4, 2]$ & $[9, 1, 1, 1, -5, 1]$ & $[0, 1, 2, 1, 1, 1]$\\
&$[1, 1, 2, 1, 2, 1]$ & $[0, 4, 1, 1, -2, 1]$ & $[0, 0, 1, 0, 1, 1]$\\
&$[2, 2, 4, 2, 4, 2]$ & $[-3, 7, 1, 1, -5, 1]$ & $[1, 1, 1, 0, 1, 0]$\\
&$[3, 3, 6, 3, 6, 3]$ & $[-6, 10, 1, 1, -8, 1]$ & $[1, 1, 2, 1, 1, 1]$\\
&$[1, 2, 2, 1, 2, 1]$ & $[3, -2, 4, 1, -2, 1]$ & $[0, 1, 1, 0, 1, 0]$\\
&$[2, 4, 4, 2, 4, 2]$ & $[3, -5, 7, 1, -5, 1]$ & $[1, 1, 1, 0, 1, 0]$\\
&$[3, 6, 6, 3, 6, 3]$ & $[3, -8, 10, 1, -8, 1]$ & $[1, 2, 2, 1, 1, 0]$\\
&$[1, 2, 3, 1, 2, 1]$ & $[3, 1, -2, 4, 1, 1]$ & $[0, 0, 1, 0, 1, 1]$\\
&$[2, 4, 6, 2, 4, 2]$ & $[3, 1, -5, 7, 1, 1]$ & $[0, 1, 2, 1, 1, 1]$\\
&$[3, 6, 9, 3, 6, 3]$ & $[3, 1, -8, 10, 1, 1]$ & $[1, 1, 2, 1, 1, 1]$\\
&$[0, 1, 1, 0, 1, 1]$ & $[6, -2, 1, 4, 1, -2]$ & $[0, 0, 1, 0, 1, 1]$\\
&$[1, 1, 1, 0, 1, 1]$ & $[0, 1, 1, 4, 1, -2]$ & $[0, 0, 1, 0, 1, 1]$\\
&$[2, 2, 2, 0, 2, 2]$ & $[-3, 1, 1, 7, 1, -5]$ & $[1, 1, 1, 0, 1, 0]$\\
&$[1, 2, 2, 1, 1, 0]$ & $[3, -2, 1, 1, 1, 4]$ & $[0, 1, 1, 0, 1, 0]$\\
&$[2, 4, 4, 2, 2, 0]$ & $[3, -5, 1, 1, 1, 7]$ & $[1, 1, 1, 0, 1, 0]$\\
&$[1, 1, 0, 0, 0, 0]$ & $[0, -2, 4, 1, 1, 1]$ & $[1, 0, 0, 0, 0, 0]$\\
&$[1, 1, 2, 1, 1, 0]$ & $[0, 4, -2, 1, 1, 4]$ & $[1, 0, 0, 0, 0, 0]$\\
&$[2, 2, 4, 2, 2, 0]$ & $[-3, 7, -5, 1, 1, 7]$ & $[1, 1, 1, 0, 1, 0]$\\
&$[1, 1, 1, 0, 0, 0]$ & $[0, 1, -2, 4, 4, 1]$ & $[1, 0, 0, 0, 0, 0]$\\
&$[1, 1, 1, 1, 1, 0]$ & $[0, 1, 4, -2, -2, 4]$ & $[1, 0, 0, 0, 0, 0]$\\
&$[2, 2, 2, 2, 2, 0]$ & $[-3, 1, 7, -5, -5, 7]$ & $[1, 1, 1, 0, 1, 0]$\\
&$[1, 1, 1, 0, 1, 0]$ & $[0, 1, 1, 4, -2, 4]$ & $[1, 0, 0, 0, 0, 0]$\\
&$[1, 1, 1, 1, 0, 0]$ & $[0, 1, 1, -2, 4, 1]$ & $[1, 0, 0, 0, 0, 0]$\\
&$[0, 1, 2, 1, 1, 0]$ & $[6, 1, -2, 1, 1, 4]$ & $[0, 1, 1, 0, 1, 0]$\\
&$[0, 1, 1, 1, 1, 0]$ & $[6, -2, 4, -2, -2, 4]$ & $[0, 1, 1, 0, 1, 0]$\\
\hline
\end{tabular}
\caption{$p=3$, relevant root multiples for $2\omega_1$}\label{small3table2}
\end{table}
\begin{table}[f]
\centering
\begin{tabular}{|c|c|c|c|c|} 
\hline
$\lambda$ & $m\alpha$ & &\\
\hline
& & $\lambda+\rho-pm\alpha$ & $\beta$ \\
\cline{2-4}
$\omega_1$
&$[1, 2, 3, 2, 2, 1]$ & $[2, 1, 1, -4, 1, 1]$ & $[0, 1, 1, 1, 1, 1]$\\
&$[2, 4, 6, 4, 4, 2]$ & $[2, 1, 1, -9, 1, 1]$ & $[1, 2, 2, 1, 2, 1]$\\
&$[1, 1, 1, 1, 1, 1]$ & $[-3, 1, 6, -4, 1, -4]$ & $[0, 1, 1, 1, 1, 1]$\\
&$[0, 1, 2, 1, 1, 1]$ & $[7, 1, -4, 1, 6, -4]$ & $[0, 1, 1, 1, 1, 1]$\\
&$[1, 1, 2, 1, 1, 1]$ & $[-3, 6, -4, 1, 6, -4]$ & $[1, 1, 1, 1, 0, 0]$\\
&$[1, 2, 2, 1, 1, 1]$ & $[2, -4, 1, 1, 6, -4]$ & $[0, 1, 1, 1, 1, 1]$\\
&$[0, 1, 2, 1, 2, 1]$ & $[7, 1, 1, 1, -4, 1]$ & $[0, 1, 1, 1, 1, 1]$\\
&$[1, 1, 2, 1, 2, 1]$ & $[-3, 6, 1, 1, -4, 1]$ & $[1, 1, 1, 0, 1, 0]$\\
&$[1, 2, 2, 1, 2, 1]$ & $[2, -4, 6, 1, -4, 1]$ & $[0, 1, 1, 1, 1, 1]$\\
&$[1, 2, 3, 1, 2, 1]$ & $[2, 1, -4, 6, 1, 1]$ & $[1, 1, 1, 0, 1, 0]$\\
&$[2, 4, 6, 2, 4, 2]$ & $[2, 1, -9, 11, 1, 1]$ & $[1, 2, 2, 1, 2, 1]$\\
&$[1, 1, 1, 0, 1, 1]$ & $[-3, 1, 1, 6, 1, -4]$ & $[1, 1, 1, 0, 1, 0]$\\
&$[1, 2, 2, 1, 1, 0]$ & $[2, -4, 1, 1, 1, 6]$ & $[1, 1, 1, 0, 1, 0]$\\
&$[1, 1, 2, 1, 1, 0]$ & $[-3, 6, -4, 1, 1, 6]$ & $[1, 1, 1, 0, 1, 0]$\\
&$[1, 1, 1, 1, 1, 0]$ & $[-3, 1, 6, -4, -4, 6]$ & $[1, 1, 1, 0, 1, 0]$\\
\hline
& & $\lambda+\rho-pm\alpha$ & $\beta$ \\
\cline{2-4}
$2\omega_1$
&$[1, 2, 3, 2, 2, 1]$ & $[3, 1, 1, -4, 1, 1]$ & $[0, 1, 1, 1, 1, 1]$\\
&$[2, 4, 6, 4, 4, 2]$ & $[3, 1, 1, -9, 1, 1]$ & $[1, 2, 2, 1, 1, 1]$\\
&$[1, 1, 1, 1, 1, 1]$ & $[-2, 1, 6, -4, 1, -4]$ & $[0, 1, 1, 1, 1, 1]$\\
&$[0, 1, 2, 1, 1, 1]$ & $[8, 1, -4, 1, 6, -4]$ & $[0, 1, 1, 1, 1, 1]$\\
&$[1, 1, 2, 1, 1, 1]$ & $[-2, 6, -4, 1, 6, -4]$ & $[1, 1, 1, 0, 0, 0]$\\
&$[1, 2, 2, 1, 1, 1]$ & $[3, -4, 1, 1, 6, -4]$ & $[0, 1, 1, 1, 1, 1]$\\
&$[0, 1, 2, 1, 2, 1]$ & $[8, 1, 1, 1, -4, 1]$ & $[0, 1, 1, 1, 1, 1]$\\
&$[1, 1, 2, 1, 2, 1]$ & $[-2, 6, 1, 1, -4, 1]$ & $[1, 1, 2, 1, 2, 1]$\\
&$[1, 2, 2, 1, 2, 1]$ & $[3, -4, 6, 1, -4, 1]$ & $[0, 1, 1, 1, 1, 1]$\\
&$[2, 4, 4, 2, 4, 2]$ & $[3, -9, 11, 1, -9, 1]$ & $[1, 2, 2, 1, 1, 1]$\\
&$[1, 2, 3, 1, 2, 1]$ & $[3, 1, -4, 6, 1, 1]$ & $[1, 1, 1, 0, 0, 0]$\\
&$[2, 4, 6, 2, 4, 2]$ & $[3, 1, -9, 11, 1, 1]$ & $[1, 2, 2, 1, 1, 1]$\\
&$[1, 1, 1, 0, 1, 1]$ & $[-2, 1, 1, 6, 1, -4]$ & $[1, 1, 1, 0, 0, 0]$\\
&$[1, 2, 2, 1, 1, 0]$ & $[3, -4, 1, 1, 1, 6]$ & $[1, 1, 1, 0, 0, 0]$\\
&$[1, 1, 2, 1, 1, 0]$ & $[-2, 6, -4, 1, 1, 6]$ & $[1, 1, 1, 0, 0, 0]$\\
&$[1, 1, 1, 0, 1, 0]$ & $[-2, 1, 1, 6, -4, 6]$ & $[1, 1, 1, 0, 0, 0]$\\
&$[1, 1, 1, 1, 0, 0]$ & $[-2, 1, 1, -4, 6, 1]$ & $[1, 1, 1, 0, 0, 0]$\\
\cline{2-4}
& & $w$ & $w(\lambda+\rho-pm\alpha)-\rho$ \\
\cline{2-4}
&$[1, 1, 1, 1, 1, 0]$ & $[1, 2, 3, 5, 4, 2, 1]$ & $[0, 0, 0, 0, 0, 1]$\\
\hline
\end{tabular}
\caption{$p=5$, relevant root multiples for $\omega_1$ and $2\omega_1$.}\label{small5table1}
\end{table}
\begin{table}[f]
\centering
\begin{tabular}{|c|c|c|c|c|} 
\hline
$\lambda$ & $m\alpha$ & &\\
\hline
& & $\lambda+\rho-pm\alpha$ & $\beta$ \\
\cline{2-4}
$3\omega_1$
&$[1, 2, 3, 2, 2, 1]$ & $[4, 1, 1, -4, 1, 1]$ & $[0, 1, 1, 1, 1, 1]$\\
&$[2, 4, 6, 4, 4, 2]$ & $[4, 1, 1, -9, 1, 1]$ & $[1, 1, 2, 1, 1, 1]$\\
&$[1, 1, 1, 1, 1, 1]$ & $[-1, 1, 6, -4, 1, -4]$ & $[0, 1, 1, 1, 1, 1]$\\
&$[0, 1, 2, 1, 1, 1]$ & $[9, 1, -4, 1, 6, -4]$ & $[0, 1, 1, 1, 1, 1]$\\
&$[1, 1, 2, 1, 1, 1]$ & $[-1, 6, -4, 1, 6, -4]$ & $[1, 1, 2, 1, 1, 1]$\\
&$[1, 2, 2, 1, 1, 1]$ & $[4, -4, 1, 1, 6, -4]$ & $[0, 1, 1, 1, 1, 1]$\\
&$[2, 4, 4, 2, 2, 2]$ & $[4, -9, 1, 1, 11, -9]$ & $[1, 1, 2, 1, 1, 1]$\\
&$[0, 1, 2, 1, 2, 1]$ & $[9, 1, 1, 1, -4, 1]$ & $[0, 1, 1, 1, 1, 1]$\\
&$[2, 2, 4, 2, 4, 2]$ & $[-6, 11, 1, 1, -9, 1]$ & $[1, 1, 2, 1, 1, 1]$\\
&$[1, 2, 2, 1, 2, 1]$ & $[4, -4, 6, 1, -4, 1]$ & $[0, 1, 1, 1, 1, 1]$\\
&$[2, 4, 4, 2, 4, 2]$ & $[4, -9, 11, 1, -9, 1]$ & $[1, 2, 2, 1, 1, 0]$\\
&$[1, 2, 3, 1, 2, 1]$ & $[4, 1, -4, 6, 1, 1]$ & $[0, 1, 2, 1, 1, 0]$\\
&$[2, 4, 6, 2, 4, 2]$ & $[4, 1, -9, 11, 1, 1]$ & $[1, 1, 2, 1, 1, 1]$\\
&$[1, 1, 1, 0, 1, 1]$ & $[-1, 1, 1, 6, 1, -4]$ & $[1, 1, 0, 0, 0, 0]$\\
&$[1, 2, 2, 1, 1, 0]$ & $[4, -4, 1, 1, 1, 6]$ & $[1, 2, 2, 1, 1, 0]$\\
&$[1, 1, 2, 1, 1, 0]$ & $[-1, 6, -4, 1, 1, 6]$ & $[0, 1, 2, 1, 1, 0]$\\
&$[1, 1, 1, 0, 0, 0]$ & $[-1, 1, -4, 6, 6, 1]$ & $[1, 1, 0, 0, 0, 0]$\\
&$[1, 1, 1, 1, 1, 0]$ & $[-1, 1, 6, -4, -4, 6]$ & $[1, 1, 0, 0, 0, 0]$\\
&$[1, 1, 1, 0, 1, 0]$ & $[-1, 1, 1, 6, -4, 6]$ & $[1, 1, 0, 0, 0, 0]$\\
&$[1, 1, 1, 1, 0, 0]$ & $[-1, 1, 1, -4, 6, 1]$ & $[1, 1, 0, 0, 0, 0]$\\
\cline{2-4}
& & $w$ & $w(\lambda+\rho-pm\alpha)-\rho$ \\
\cline{2-4}
&$[1, 1, 2, 1, 2, 1]$ & $[3, 5, 4, 3, 6, 5, 1]$ & $[0, 0, 0, 0, 0, 0]$\\
\hline
& & $\lambda+\rho-pm\alpha$ & $\beta$ \\
\cline{2-4}
$4\omega_1$
&$[1, 2, 3, 2, 2, 1]$ & $[5, 1, 1, -4, 1, 1]$ & $[0, 1, 1, 1, 1, 1]$\\
&$[2, 4, 6, 4, 4, 2]$ & $[5, 1, 1, -9, 1, 1]$ & $[1, 1, 1, 1, 1, 1]$\\
&$[1, 1, 1, 1, 1, 1]$ & $[0, 1, 6, -4, 1, -4]$ & $[0, 1, 1, 1, 1, 1]$\\
&$[0, 1, 2, 1, 1, 1]$ & $[10, 1, -4, 1, 6, -4]$ & $[0, 1, 1, 1, 1, 1]$\\
&$[1, 1, 2, 1, 1, 1]$ & $[0, 6, -4, 1, 6, -4]$ & $[1, 0, 0, 0, 0, 0]$\\
&$[2, 2, 4, 2, 2, 2]$ & $[-5, 11, -9, 1, 11, -9]$ & $[1, 1, 1, 1, 1, 1]$\\
&$[1, 2, 2, 1, 1, 1]$ & $[5, -4, 1, 1, 6, -4]$ & $[0, 1, 1, 1, 1, 1]$\\
&$[2, 4, 4, 2, 2, 2]$ & $[5, -9, 1, 1, 11, -9]$ & $[1, 1, 1, 1, 1, 1]$\\
&$[0, 1, 2, 1, 2, 1]$ & $[10, 1, 1, 1, -4, 1]$ & $[0, 1, 1, 1, 1, 1]$\\
&$[1, 1, 2, 1, 2, 1]$ & $[0, 6, 1, 1, -4, 1]$ & $[1, 0, 0, 0, 0, 0]$\\
&$[2, 2, 4, 2, 4, 2]$ & $[-5, 11, 1, 1, -9, 1]$ & $[1, 1, 1, 1, 1, 1]$\\
&$[1, 2, 2, 1, 2, 1]$ & $[5, -4, 6, 1, -4, 1]$ & $[0, 1, 1, 1, 1, 1]$\\
&$[2, 4, 4, 2, 4, 2]$ & $[5, -9, 11, 1, -9, 1]$ & $[1, 1, 1, 1, 1, 1]$\\
&$[1, 2, 3, 1, 2, 1]$ & $[5, 1, -4, 6, 1, 1]$ & $[0, 1, 2, 1, 1, 0]$\\
&$[2, 4, 6, 2, 4, 2]$ & $[5, 1, -9, 11, 1, 1]$ & $[1, 1, 2, 1, 1, 0]$\\
&$[1, 1, 1, 0, 1, 1]$ & $[0, 1, 1, 6, 1, -4]$ & $[1, 0, 0, 0, 0, 0]$\\
&$[1, 2, 2, 1, 1, 0]$ & $[5, -4, 1, 1, 1, 6]$ & $[0, 1, 2, 1, 1, 0]$\\
&$[2, 4, 4, 2, 2, 0]$ & $[5, -9, 1, 1, 1, 11]$ & $[1, 1, 2, 1, 1, 0]$\\
&$[1, 1, 0, 0, 0, 0]$ & $[0, -4, 6, 1, 1, 1]$ & $[1, 0, 0, 0, 0, 0]$\\
&$[1, 1, 2, 1, 1, 0]$ & $[0, 6, -4, 1, 1, 6]$ & $[1, 0, 0, 0, 0, 0]$\\
&$[1, 1, 1, 0, 0, 0]$ & $[0, 1, -4, 6, 6, 1]$ & $[1, 0, 0, 0, 0, 0]$\\
&$[1, 1, 1, 1, 1, 0]$ & $[0, 1, 6, -4, -4, 6]$ & $[1, 0, 0, 0, 0, 0]$\\
&$[1, 1, 1, 0, 1, 0]$ & $[0, 1, 1, 6, -4, 6]$ & $[1, 0, 0, 0, 0, 0]$\\
&$[1, 1, 1, 1, 0, 0]$ & $[0, 1, 1, -4, 6, 1]$ & $[1, 0, 0, 0, 0, 0]$\\
\hline
\end{tabular}
\caption{$p=5$, relevant root multiples for $3\omega_1$ and $4\omega_1$.}\label{small5table2}
\end{table}
\begin{table}[f]
\centering
\begin{tabular}{|c|c|c|c|c|} 
\hline
$\lambda$ & $m\alpha$ & &\\
\hline
& & $\lambda+\rho-pm\alpha$ & $\beta$ \\
\cline{2-4}
$\omega_1$
&$[1, 2, 3, 2, 2, 1]$ & $[2, 1, 1, -6, 1, 1]$ & $[1, 1, 1, 1, 1, 1]$\\
&$[1, 1, 2, 1, 1, 1]$ & $[-5, 8, -6, 1, 8, -6]$ & $[1, 1, 1, 1, 1, 1]$\\
&$[1, 2, 2, 1, 1, 1]$ & $[2, -6, 1, 1, 8, -6]$ & $[1, 1, 1, 1, 1, 1]$\\
&$[1, 1, 2, 1, 2, 1]$ & $[-5, 8, 1, 1, -6, 1]$ & $[1, 1, 1, 1, 1, 1]$\\
&$[1, 2, 2, 1, 2, 1]$ & $[2, -6, 8, 1, -6, 1]$ & $[1, 1, 1, 1, 1, 1]$\\
&$[1, 2, 3, 1, 2, 1]$ & $[2, 1, -6, 8, 1, 1]$ & $[0, 1, 2, 1, 2, 1]$\\
&$[1, 2, 2, 1, 1, 0]$ & $[2, -6, 1, 1, 1, 8]$ & $[1, 1, 2, 1, 1, 0]$\\
\hline
& & $\lambda+\rho-pm\alpha$ & $\beta$ \\
\cline{2-4}
$2\omega_1$
&$[1, 2, 3, 2, 2, 1]$ & $[3, 1, 1, -6, 1, 1]$ & $[0, 1, 2, 1, 2, 1]$\\
&$[1, 1, 1, 1, 1, 1]$ & $[-4, 1, 8, -6, 1, -6]$ & $[1, 1, 1, 0, 1, 1]$\\
&$[1, 1, 2, 1, 1, 1]$ & $[-4, 8, -6, 1, 8, -6]$ & $[1, 1, 1, 0, 1, 1]$\\
&$[1, 2, 2, 1, 1, 1]$ & $[3, -6, 1, 1, 8, -6]$ & $[1, 1, 1, 0, 1, 1]$\\
&$[1, 1, 2, 1, 2, 1]$ & $[-4, 8, 1, 1, -6, 1]$ & $[0, 1, 2, 1, 2, 1]$\\
&$[1, 2, 2, 1, 2, 1]$ & $[3, -6, 8, 1, -6, 1]$ & $[0, 1, 2, 1, 2, 1]$\\
&$[1, 2, 3, 1, 2, 1]$ & $[3, 1, -6, 8, 1, 1]$ & $[0, 1, 2, 1, 2, 1]$\\
&$[1, 2, 2, 1, 1, 0]$ & $[3, -6, 1, 1, 1, 8]$ & $[1, 1, 1, 1, 1, 0]$\\
&$[1, 1, 2, 1, 1, 0]$ & $[-4, 8, -6, 1, 1, 8]$ & $[1, 1, 1, 1, 1, 0]$\\
\hline
& & $\lambda+\rho-pm\alpha$ & $\beta$ \\
\cline{2-4}
$3\omega_1$
&$[1, 2, 3, 2, 2, 1]$ & $[4, 1, 1, -6, 1, 1]$ & $[1, 2, 3, 2, 2, 1]$\\
&$[1, 1, 1, 1, 1, 1]$ & $[-3, 1, 8, -6, 1, -6]$ & $[1, 1, 1, 1, 0, 0]$\\
&$[1, 1, 2, 1, 1, 1]$ & $[-3, 8, -6, 1, 8, -6]$ & $[1, 1, 1, 1, 0, 0]$\\
&$[1, 2, 2, 1, 1, 1]$ & $[4, -6, 1, 1, 8, -6]$ & $[1, 1, 1, 1, 0, 0]$\\
&$[1, 1, 2, 1, 2, 1]$ & $[-3, 8, 1, 1, -6, 1]$ & $[0, 1, 2, 1, 2, 1]$\\
&$[1, 2, 2, 1, 2, 1]$ & $[4, -6, 8, 1, -6, 1]$ & $[0, 1, 2, 1, 2, 1]$\\
&$[1, 2, 3, 1, 2, 1]$ & $[4, 1, -6, 8, 1, 1]$ & $[0, 1, 2, 1, 2, 1]$\\
&$[1, 1, 1, 0, 1, 1]$ & $[-3, 1, 1, 8, 1, -6]$ & $[1, 1, 1, 0, 1, 0]$\\
&$[1, 2, 2, 1, 1, 0]$ & $[4, -6, 1, 1, 1, 8]$ & $[1, 1, 1, 0, 1, 0]$\\
&$[1, 1, 2, 1, 1, 0]$ & $[-3, 8, -6, 1, 1, 8]$ & $[1, 1, 1, 0, 1, 0]$\\
&$[1, 1, 1, 1, 1, 0]$ & $[-3, 1, 8, -6, -6, 8]$ & $[1, 1, 1, 0, 1, 0]$\\
\hline
\end{tabular}
\caption{$p=7$, relevant root multiples for $r\omega_1$, $1\leq r\leq 3$.}\label{small7table1}
\end{table}
\begin{table}[f]
\centering
\begin{tabular}{|c|c|c|c|c|} 
\hline
$\lambda$ & $m\alpha$ & &\\
\hline
& & $\lambda+\rho-pm\alpha$ & $\beta$ \\
\cline{2-4}
$4\omega_1$
&$[1, 2, 3, 2, 2, 1]$ & $[5, 1, 1, -6, 1, 1]$ & $[0, 1, 2, 1, 2, 1]$\\
&$[2, 4, 6, 4, 4, 2]$ & $[5, 1, 1, -13, 1, 1]$ & $[1, 2, 3, 1, 2, 1]$\\
&$[1, 1, 2, 1, 1, 1]$ & $[-2, 8, -6, 1, 8, -6]$ & $[1, 1, 1, 0, 0, 0]$\\
&$[1, 2, 2, 1, 1, 1]$ & $[5, -6, 1, 1, 8, -6]$ & $[1, 1, 1, 0, 0, 0]$\\
&$[1, 1, 2, 1, 2, 1]$ & $[-2, 8, 1, 1, -6, 1]$ & $[0, 1, 2, 1, 2, 1]$\\
&$[1, 2, 2, 1, 2, 1]$ & $[5, -6, 8, 1, -6, 1]$ & $[0, 1, 2, 1, 2, 1]$\\
&$[1, 2, 3, 1, 2, 1]$ & $[5, 1, -6, 8, 1, 1]$ & $[0, 1, 2, 1, 2, 1]$\\
&$[1, 1, 1, 0, 1, 1]$ & $[-2, 1, 1, 8, 1, -6]$ & $[1, 1, 1, 0, 0, 0]$\\
&$[1, 2, 2, 1, 1, 0]$ & $[5, -6, 1, 1, 1, 8]$ & $[1, 1, 1, 0, 0, 0]$\\
&$[1, 1, 2, 1, 1, 0]$ & $[-2, 8, -6, 1, 1, 8]$ & $[1, 1, 1, 0, 0, 0]$\\
&$[1, 1, 1, 0, 1, 0]$ & $[-2, 1, 1, 8, -6, 8]$ & $[1, 1, 1, 0, 0, 0]$\\
&$[1, 1, 1, 1, 0, 0]$ & $[-2, 1, 1, -6, 8, 1]$ & $[1, 1, 1, 0, 0, 0]$\\
\cline{2-4}
& & $w$ & $w(\lambda+\rho-pm\alpha)-\rho$ \\
\cline{2-4}
&$[1, 1, 1, 1, 1, 1]$ & $[1, 2, 3, 5, 6, 4, 2, 1]$ & $[1, 0, 0, 1, 0, 0]$\\
&$[1, 1, 1, 1, 1, 0]$ & $[1, 2, 3, 5, 4, 2, 1]$ & $[2, 0, 0, 0, 0, 1]$\\
\hline
& & $\lambda+\rho-pm\alpha$ & $\beta$ \\
\cline{2-4}
$5\omega_1$
&$[1, 2, 3, 2, 2, 1]$ & $[6, 1, 1, -6, 1, 1]$ & $[0, 1, 2, 1, 2, 1]$\\
&$[2, 4, 6, 4, 4, 2]$ & $[6, 1, 1, -13, 1, 1]$ & $[1, 2, 2, 1, 2, 1]$\\
&$[1, 1, 1, 1, 1, 1]$ & $[-1, 1, 8, -6, 1, -6]$ & $[1, 1, 0, 0, 0, 0]$\\
&$[1, 2, 2, 1, 1, 1]$ & $[6, -6, 1, 1, 8, -6]$ & $[1, 1, 0, 0, 0, 0]$\\
&$[1, 1, 2, 1, 2, 1]$ & $[-1, 8, 1, 1, -6, 1]$ & $[0, 1, 2, 1, 2, 1]$\\
&$[1, 2, 2, 1, 2, 1]$ & $[6, -6, 8, 1, -6, 1]$ & $[0, 1, 2, 1, 2, 1]$\\
&$[1, 2, 3, 1, 2, 1]$ & $[6, 1, -6, 8, 1, 1]$ & $[0, 1, 2, 1, 2, 1]$\\
&$[2, 4, 6, 2, 4, 2]$ & $[6, 1, -13, 15, 1, 1]$ & $[1, 2, 2, 1, 2, 1]$\\
&$[1, 1, 1, 0, 1, 1]$ & $[-1, 1, 1, 8, 1, -6]$ & $[1, 1, 0, 0, 0, 0]$\\
&$[1, 2, 2, 1, 1, 0]$ & $[6, -6, 1, 1, 1, 8]$ & $[1, 1, 0, 0, 0, 0]$\\
&$[1, 1, 1, 0, 0, 0]$ & $[-1, 1, -6, 8, 8, 1]$ & $[1, 1, 0, 0, 0, 0]$\\
&$[1, 1, 1, 1, 1, 0]$ & $[-1, 1, 8, -6, -6, 8]$ & $[1, 1, 0, 0, 0, 0]$\\
&$[1, 1, 1, 0, 1, 0]$ & $[-1, 1, 1, 8, -6, 8]$ & $[1, 1, 0, 0, 0, 0]$\\
&$[1, 1, 1, 1, 0, 0]$ & $[-1, 1, 1, -6, 8, 1]$ & $[1, 1, 0, 0, 0, 0]$\\
\cline{2-4}
& & $w$ & $w(\lambda+\rho-pm\alpha)-\rho$ \\
\cline{2-4}
&$[1, 1, 2, 1, 1, 1]$ & $[1, 2, 3, 5, 6, 4, 3, 1]$ & $[0, 0, 0, 1, 0, 1]$\\
&$[1, 1, 2, 1, 1, 0]$ & $[1, 2, 3, 5, 4, 3, 1]$ & $[1, 0, 0, 0, 0, 2]$\\
\hline
\end{tabular}
\caption{$p=7$, relevant root multiples for $r\omega_1$, $4\leq r\leq 5$.}\label{small7table2}
\end{table}
\begin{table}[f]
\centering
\begin{tabular}{|c|c|c|c|c|} 
\hline
$\lambda$ & $m\alpha$ & &\\
\cline{2-4}
& & $\lambda+\rho-pm\alpha$ & $\beta$ \\
\cline{2-4}
$6\omega_1$
&$[1, 2, 3, 2, 2, 1]$ & $[7, 1, 1, -6, 1, 1]$ & $[0, 1, 2, 1, 2, 1]$\\
&$[2, 4, 6, 4, 4, 2]$ & $[7, 1, 1, -13, 1, 1]$ & $[1, 2, 2, 1, 1, 1]$\\
&$[1, 1, 1, 1, 1, 1]$ & $[0, 1, 8, -6, 1, -6]$ & $[1, 0, 0, 0, 0, 0]$\\
&$[1, 1, 2, 1, 1, 1]$ & $[0, 8, -6, 1, 8, -6]$ & $[1, 0, 0, 0, 0, 0]$\\
&$[1, 2, 2, 1, 1, 1]$ & $[7, -6, 1, 1, 8, -6]$ & $[1, 2, 2, 1, 1, 1]$\\
&$[1, 1, 2, 1, 2, 1]$ & $[0, 8, 1, 1, -6, 1]$ & $[0, 1, 2, 1, 2, 1]$\\
&$[1, 2, 2, 1, 2, 1]$ & $[7, -6, 8, 1, -6, 1]$ & $[0, 1, 2, 1, 2, 1]$\\
&$[2, 4, 4, 2, 4, 2]$ & $[7, -13, 15, 1, -13, 1]$ & $[1, 2, 2, 1, 1, 1]$\\
&$[1, 2, 3, 1, 2, 1]$ & $[7, 1, -6, 8, 1, 1]$ & $[0, 1, 2, 1, 2, 1]$\\
&$[2, 4, 6, 2, 4, 2]$ & $[7, 1, -13, 15, 1, 1]$ & $[1, 2, 2, 1, 1, 1]$\\
&$[1, 1, 1, 0, 1, 1]$ & $[0, 1, 1, 8, 1, -6]$ & $[1, 0, 0, 0, 0, 0]$\\
&$[1, 1, 0, 0, 0, 0]$ & $[0, -6, 8, 1, 1, 1]$ & $[1, 0, 0, 0, 0, 0]$\\
&$[1, 1, 2, 1, 1, 0]$ & $[0, 8, -6, 1, 1, 8]$ & $[1, 0, 0, 0, 0, 0]$\\
&$[1, 1, 1, 0, 0, 0]$ & $[0, 1, -6, 8, 8, 1]$ & $[1, 0, 0, 0, 0, 0]$\\
&$[1, 1, 1, 1, 1, 0]$ & $[0, 1, 8, -6, -6, 8]$ & $[1, 0, 0, 0, 0, 0]$\\
&$[1, 1, 1, 0, 1, 0]$ & $[0, 1, 1, 8, -6, 8]$ & $[1, 0, 0, 0, 0, 0]$\\
&$[1, 1, 1, 1, 0, 0]$ & $[0, 1, 1, -6, 8, 1]$ & $[1, 0, 0, 0, 0, 0]$\\
\cline{2-4}
& & $w$ & $w(\lambda+\rho-pm\alpha)-\rho$ \\
\cline{2-4}
&$[1, 2, 2, 1, 1, 0]$ & $[1, 2, 3, 5, 4, 3, 2]$ & $[0, 0, 0, 0, 0, 3]$\\
\hline
\end{tabular}
\caption{$p=7$, relevant root multiples for $6\omega_1$.}\label{small7table3}
\end{table}
\begin{table}[f]
\centering
\begin{tabular}{|c|c|c|c|c|} 
\hline
$\lambda$ & $m\alpha$ & &\\
\hline
& & $\lambda+\rho-pm\alpha$ & $\beta$ \\
\cline{2-4}
$2\omega_1+\omega_6$
&$[1, 2, 3, 2, 2, 1]$ & $[3, 1, 1, -6, 1, 2]$ & $[1, 2, 3, 2, 2, 1]$\\
&$[1, 1, 1, 1, 1, 1]$ & $[-4, 1, 8, -6, 1, -5]$ & $[1, 1, 1, 1, 1, 0]$\\
&$[1, 1, 2, 1, 1, 1]$ & $[-4, 8, -6, 1, 8, -5]$ & $[0, 1, 2, 1, 1, 1]$\\
&$[1, 2, 2, 1, 1, 1]$ & $[3, -6, 1, 1, 8, -5]$ & $[0, 1, 2, 1, 1, 1]$\\
&$[0, 1, 2, 1, 2, 1]$ & $[10, 1, 1, 1, -6, 2]$ & $[0, 1, 2, 1, 1, 1]$\\
&$[1, 1, 2, 1, 2, 1]$ & $[-4, 8, 1, 1, -6, 2]$ & $[1, 1, 1, 1, 1, 0]$\\
&$[1, 2, 2, 1, 2, 1]$ & $[3, -6, 8, 1, -6, 2]$ & $[1, 1, 1, 1, 1, 0]$\\
&$[1, 2, 3, 1, 2, 1]$ & $[3, 1, -6, 8, 1, 2]$ & $[0, 1, 2, 1, 1, 1]$\\
&$[1, 2, 2, 1, 1, 0]$ & $[3, -6, 1, 1, 1, 9]$ & $[1, 1, 1, 1, 1, 0]$\\
&$[1, 1, 2, 1, 1, 0]$ & $[-4, 8, -6, 1, 1, 9]$ & $[1, 1, 1, 1, 1, 0]$\\
\cline{2-4}
& & $w$ & $w(\lambda+\rho-pm\alpha)-\rho$ \\
\cline{2-4}
&$[1, 1, 1, 0, 1, 1]$ & $[1, 2, 3, 5, 6, 5, 3, 2, 1]$ & $[1, 0, 0, 1, 0, 0]$\\
\hline
& & $\lambda+\rho-pm\alpha$ & $\beta$ \\
\cline{2-4}
$\omega_1+\omega_4$
&$[1, 2, 3, 2, 2, 1]$ & $[2, 1, 1, -5, 1, 1]$ & $[1, 2, 3, 2, 2, 1]$\\
&$[1, 1, 1, 1, 1, 1]$ & $[-5, 1, 8, -5, 1, -6]$ & $[1, 1, 1, 1, 1, 0]$\\
&$[1, 1, 2, 1, 1, 1]$ & $[-5, 8, -6, 2, 8, -6]$ & $[0, 1, 2, 1, 1, 1]$\\
&$[1, 2, 2, 1, 1, 1]$ & $[2, -6, 1, 2, 8, -6]$ & $[0, 1, 2, 1, 1, 1]$\\
&$[0, 1, 2, 1, 2, 1]$ & $[9, 1, 1, 2, -6, 1]$ & $[0, 1, 2, 1, 1, 1]$\\
&$[1, 1, 2, 1, 2, 1]$ & $[-5, 8, 1, 2, -6, 1]$ & $[1, 1, 1, 1, 1, 0]$\\
&$[1, 2, 2, 1, 2, 1]$ & $[2, -6, 8, 2, -6, 1]$ & $[1, 1, 1, 1, 1, 0]$\\
&$[1, 2, 3, 1, 2, 1]$ & $[2, 1, -6, 9, 1, 1]$ & $[0, 1, 2, 1, 1, 1]$\\
&$[1, 2, 2, 1, 1, 0]$ & $[2, -6, 1, 2, 1, 8]$ & $[1, 1, 1, 1, 1, 0]$\\
&$[1, 1, 2, 1, 1, 0]$ & $[-5, 8, -6, 2, 1, 8]$ & $[1, 1, 1, 1, 1, 0]$\\
\hline
\end{tabular}
\caption{$p=7$, second iteration relevant root multiples for $4\omega_1$.}\label{small7table4}
\end{table}
\begin{table}[f]
\centering
\begin{tabular}{|c|c|c|c|c|} 
\hline
$\lambda$ & $m\alpha$ & &\\
\hline
& & $\lambda+\rho-pm\alpha$ & $\beta$ \\
\cline{2-4}
$\omega_1+2\omega_6$
&$[1, 2, 3, 2, 2, 1]$ & $[2, 1, 1, -6, 1, 3]$ & $[1, 2, 3, 2, 2, 1]$\\
&$[1, 1, 1, 1, 1, 1]$ & $[-5, 1, 8, -6, 1, -4]$ & $[0, 1, 1, 1, 1, 1]$\\
&$[0, 1, 2, 1, 1, 1]$ & $[9, 1, -6, 1, 8, -4]$ & $[0, 1, 1, 1, 1, 1]$\\
&$[1, 1, 2, 1, 1, 1]$ & $[-5, 8, -6, 1, 8, -4]$ & $[1, 1, 2, 1, 1, 0]$\\
&$[1, 2, 2, 1, 1, 1]$ & $[2, -6, 1, 1, 8, -4]$ & $[0, 1, 1, 1, 1, 1]$\\
&$[0, 1, 2, 1, 2, 1]$ & $[9, 1, 1, 1, -6, 3]$ & $[0, 1, 1, 1, 1, 1]$\\
&$[1, 1, 2, 1, 2, 1]$ & $[-5, 8, 1, 1, -6, 3]$ & $[1, 1, 2, 1, 1, 0]$\\
&$[1, 2, 2, 1, 2, 1]$ & $[2, -6, 8, 1, -6, 3]$ & $[0, 1, 1, 1, 1, 1]$\\
&$[1, 2, 3, 1, 2, 1]$ & $[2, 1, -6, 8, 1, 3]$ & $[1, 1, 2, 1, 1, 0]$\\
&$[1, 2, 2, 1, 1, 0]$ & $[2, -6, 1, 1, 1, 10]$ & $[1, 1, 2, 1, 1, 0]$\\
\cline{2-4}
& & $w$ & $w(\lambda+\rho-pm\alpha)-\rho$ \\
\cline{2-4}
&$[1, 1, 1, 0, 1, 1]$ & $[1, 2, 3, 5, 6, 5, 3, 2, 1]$ & $[0, 0, 0, 1, 0, 1]$\\
\hline
& & $\lambda+\rho-pm\alpha$ & $\beta$ \\
\cline{2-4}
$\omega_4+\omega_6$
&$[1, 2, 3, 2, 2, 1]$ & $[1, 1, 1, -5, 1, 2]$ & $[1, 2, 3, 2, 2, 1]$\\
&$[1, 1, 1, 1, 1, 1]$ & $[-6, 1, 8, -5, 1, -5]$ & $[0, 1, 1, 1, 1, 1]$\\
&$[0, 1, 2, 1, 1, 1]$ & $[8, 1, -6, 2, 8, -5]$ & $[0, 1, 1, 1, 1, 1]$\\
&$[1, 1, 2, 1, 1, 1]$ & $[-6, 8, -6, 2, 8, -5]$ & $[1, 1, 2, 1, 1, 0]$\\
&$[1, 2, 2, 1, 1, 1]$ & $[1, -6, 1, 2, 8, -5]$ & $[0, 1, 1, 1, 1, 1]$\\
&$[0, 1, 2, 1, 2, 1]$ & $[8, 1, 1, 2, -6, 2]$ & $[0, 1, 1, 1, 1, 1]$\\
&$[1, 1, 2, 1, 2, 1]$ & $[-6, 8, 1, 2, -6, 2]$ & $[1, 1, 2, 1, 1, 0]$\\
&$[1, 2, 2, 1, 2, 1]$ & $[1, -6, 8, 2, -6, 2]$ & $[0, 1, 1, 1, 1, 1]$\\
&$[1, 2, 3, 1, 2, 1]$ & $[1, 1, -6, 9, 1, 2]$ & $[1, 1, 2, 1, 1, 0]$\\
&$[1, 2, 2, 1, 1, 0]$ & $[1, -6, 1, 2, 1, 9]$ & $[1, 1, 2, 1, 1, 0]$\\
\hline
\end{tabular}
\caption{$p=7$, second iteration relevant root multiples for $5\omega_1$.}\label{small7table5}
\end{table}

\subsection*{Acknowledgement} I wish to thank the referees for their
careful reading and helpful suggestions, which have improved the exposition.


\begin{thebibliography}{99}
\bibitem{Ars} O. Arslan and P. Sin, Some simple modules for classical groups and p-ranks of orthogonal and Hermitian geometries, {\it Journal of Algebra} {\bf 327} (2011) 141--169 .

\bibitem{Jantzen} J. C. Jantzen,  {\it Representations of Algebraic Groups}, Academic Press, London, 1987.

\bibitem{RR} S. Ramanan and  A. Ramanathan, Projective normality of flag varieties and Schubert varieties, {\it Inventiones Mathematicae} {\bf 79} (1985) 217--224.

\bibitem{opp} P. Sin, Oppositeness in buildings and simple modules for finite groups of Lie type, in  {\it Buildings, Finite Geometries and Groups}, Springer Proceedings in Mathematics Volume 10, (2011), 273--286.

\bibitem{St} R. Steinberg, Representations of algebraic groups, Nagoya Math. J.
{\bf 22} (1963), 33--56.

\bibitem{St2} R. Steinberg, Endomorphisms of linear algebraic groups,
Memoirs Amer. Math. Soc {\bf 80} (1968).
\end{thebibliography}
\end{document}